\documentclass[11pt,oneside,leqno]{amsart}
\usepackage{amsxtra}
\usepackage{amsopn}
\usepackage{color}
\usepackage{amsmath,amsthm,amssymb}
\usepackage{amscd}
\usepackage{amsfonts}
\usepackage{latexsym}
\usepackage{verbatim}
\usepackage{multirow}

\theoremstyle{plain}
\newtheorem{theorem}{Theorem}[section]
\newtheorem{definition}[theorem]{Definition}
\newtheorem{lemma}[theorem]{Lemma}
\newtheorem{proposition}[theorem]{Proposition}
\newtheorem{corollary}[theorem]{Corollary}
\newtheorem{remark}[theorem]{Remark}

\newtheorem{example}[theorem]{Example}

\newtheorem{remark-question}[section]{Remark-Question}

\newcommand\C{{\mathbb C}}

\newcommand\fra{{\mathfrak a}} 
\newcommand\frg{{\mathfrak g}}
\newcommand\frh{{\mathfrak h}}

\newcommand\frb{{\mathfrak b}}

\newcommand\Real{{\mathfrak R}{\frak e}\,} 
\newcommand\Imag{{\mathfrak I}{\frak m}\,}
\newcommand\nilm{\Gamma\backslash G}

\newcommand\db{{\bar{\partial}}}

\newcommand{\del}{\partial}
\newcommand{\delbar}{\overline{\del}}

\DeclareMathOperator{\imm}{im}

\begin{document}
\title[]{On the stability of compact pseudo-K\"ahler \\and neutral Calabi-Yau manifolds}

\keywords{Complex manifold; symplectic form; pseudo-K\"ahler metric; neutral Calabi-Yau structure; holomorphic deformation; nilmanifold; solvmanifold}
\subjclass[2000]{}

\author{Adela Latorre}
\address[A. Latorre]{Departamento de Matem\'atica Aplicada,
Universidad Polit\'ecnica de Madrid,
C/ Jos\'e Antonio Novais 10,
28040 Madrid, Spain}
\email{adela.latorre@upm.es}

\author{Luis Ugarte}
\address[L. Ugarte]{Departamento de Matem\'aticas\,-\,I.U.M.A.\\
Universidad de Zaragoza\\
Campus Plaza San Francisco\\
50009 Zaragoza, Spain}
\email{ugarte@unizar.es}


\maketitle

\begin{abstract}
We study the stability of compact pseudo-K\"ahler manifolds, i.e. compact complex manifolds $X$ endowed with a symplectic form compatible with the complex structure of~$X$. When the corresponding metric is positive-definite, $X$ is K\"ahler and any sufficiently small deformation of $X$ admits a K\"ahler metric by a well-known result of Kodaira and Spencer. We prove that compact pseudo-K\"ahler surfaces are also stable, 
but we show that stability fails in every complex dimension $n\geq 3$. Similar results are obtained for compact neutral K\"ahler and neutral Calabi-Yau manifolds. Finally, motivated by a question of Streets and Tian in the positive-definite case, we construct compact complex manifolds with pseudo-Hermitian-symplectic structures that
do not admit any pseudo-K\"ahler metric.
\end{abstract}

\maketitle


\setcounter{tocdepth}{3} \tableofcontents

\section{Introduction}\label{intro}

\noindent 
Let~$M$ be an even-dimensional manifold endowed with a complex structure $J$ and a symplectic
form~$F$. When~$J$ and~$F$ are compatible, i.e. the symplectic form $F$ is $J$-invariant, and the associated metric~$g$ is Riemannian, then
the manifold $(M,J,F)$ is K\"ahler. In the compact K\"ahler case, the positive-definiteness of~$g$
imposes strong topological conditions 
on the manifold~$M$; for instance, its Betti numbers
$b_{2k+1}(M)$ are even and the manifold is formal \cite{DGMS}. 
Since
there are compact complex manifolds 
with no K\"ahler metrics, many efforts have been done in understanding the properties of manifolds endowed with a pair $(J,F)$ satisfying weaker conditions than those of a K\"ahler structure. 
On the one hand, if we drop the closedness condition for the 2-form $F$, other special Hermitian structures arise, such as the 
strong K\"ahler with torsion (SKT) or balanced Hermitian structures. On the other hand, if we no longer require the positive definiteness of the metric $g$ but preserve the compatibility of the symplectic form $F$ with $J$,
then a structure called \emph{pseudo-K\"ahler} is obtained.

In this paper we focus on compact manifolds~$M$ endowed with a 
pseudo-K\"ahler structure, i.e. a pair $(J,F)$ where~$J$ is a complex structure
and~$F$ is a non-degenerate closed $2$-form such that $F(JU,JV)=F(U,V)$,
for any smooth vector fields~$U,V$ on~$M$. This is equivalent to $J$ being parallel, i.e. $\nabla J=0$, where~$\nabla$ is the Levi-Civita connection of the pseudo-Riemmanian metric $g(U,V)=F(JU,V)$ (see \cite{AFGM}).
If the real dimension of~$M$ is~$2n$, then the pseudo-K\"ahler metric has signature $(2k,2n-2k)$, where~$k=n$ corresponds to the K\"ahler case.  
In other words, we have a compact complex manifold $X=(M,J)$ with a symplectic form~$F$ of bidegree $(1,1)$ with respect to the complex structure~$J$. 
Notice that there are many compact pseudo-K\"ahler manifolds not admitting any K\"ahler metric, the simplest example being the compact complex surface known as the Kodaira-Thurston manifold~\cite{Thu}. 

Pseudo-K\"ahler structures appear in relation to other interesting structures on manifolds. 
For instance, compact complex homogeneous manifolds endowed with a pseudo-K\"ahler structure are classified in \cite{DGuan, Guan}.
Pseudo-K\"ahler Einstein metrics on compact complex surfaces 
are studied by Petean \cite{Petean}. 
More recently, it is proved in~\cite{KL} that there is a natural pseudo-K\"ahler structure on the
universal intermediate $G_2$-Jacobian $\mathcal{J}$
of the moduli space of torsion-free $G_2$-structures
on a fixed compact $7$-manifold. 
Furthermore, any $4n$-dimensional hypersymplectic manifold  
in particular has
a neutral Calabi-Yau structure
\cite{Hitchin} (see also \cite{DS}), which is a special type of pseudo-K\"ahler structure
whose underlying metric is Ricci-flat and has signature $(2n,2n)$.

Although several aspects in pseudo-K\"ahler geometry have been investigated, the stability of these structures under small holomorphic deformations of the complex manifold is, to our knowlegde, only known in the positive-definite case. Indeed, if $X$ is a compact K\"ahler manifold, then any sufficiently small deformation of $X$ admits a K\"ahler metric due to a well-known result by Kodaira and Spencer~\cite{KS60}. In this paper we focus on the stability properties of compact 
pseudo-K\"ahler manifolds, as well as manifolds with related neutral (K\"ahler and Calabi-Yau) metrics. 

We next explain in more detail the contents of the paper. 

In Section~\ref{general-results} we firstly construct compact pseudo-K\"ahler manifolds of complex dimension $n\geq 3$ that are not stable under small holomorphic deformations of the
complex structure. This result motivates the study of sufficient conditions 
for the stability of the pseudo-K\"ahler property. 
Our stability result involves the Bott-Chern cohomology of complex manifolds. We prove that if~$X_0$ is pseudo-K\"ahler and the upper-semi-continuous function $t\mapsto h^{1,1}_{\rm BC}(X_t)$ is constant, then the compact complex manifold $X_t$ admits a pseudo-K\"ahler metric for any small enough $t$ (see Proposition~\ref{BC-constant}). Combining this with a result of Teleman~\cite{Tel06} on the complex invariant 
$\Delta^2$
introduced by Angella and Tomassini in~\cite{AT4},
we prove in Theorem~\ref{stability-p-K-surfaces} that
compact pseudo-K\"ahler surfaces are stable. 

The results of Section~\ref{general-results} are illustrated with explicit constructions of pseudo-K\"ahler nilmanifolds and solvmanifolds, together with their small holomorphic deformations. In particular, we prove in Proposition~\ref{iwasawa} that the Iwasawa manifold and its small deformations do not admit any pseudo-K\"ahler metric. In contrast, 
the holomorphically parallelizable Nakamura manifold $X$ is pseudo-K\"ahler, as proved by Yamada in~\cite{Yamada-1}, and there exists a small deformation $X_{\mathbf{t}}$ of $X$ admitting pseudo-K\"ahler metrics for every $\mathbf{t}$
(see Proposition~\ref{nakamura}).

Section~\ref{neutralCY} is devoted to a special class
of pseudo-K\"ahler manifolds, namely, neutral Calabi-Yau manifolds. 
They are neutral K\"ahler manifolds, i.e. manifolds with even complex dimension $2m$ and a metric $g$ 
of signature $(2m,2m)$, that additionally have 
a nowhere vanishing form $\Phi$ of bidegree $(2m,0)$ 
satisfying $\nabla \Phi=0$, where $\nabla$ is the Levi-Civita connection of $g$. Observe that neutral Calabi-Yau manifolds are Ricci-flat. 
In Section~\ref{neutralCY-surfaces} we prove the stability of neutral K\"ahler and neutral Calabi-Yau structures on compact complex surfaces. 
However, 
our Theorem~\ref{deform-pK} shows 
that such structures are not stable on compact complex manifolds of  
any complex dimension $n\geq 4$. 
This result is in deep contrast to the case of K\"ahler Calabi-Yau manifolds, whose deformation space is unobstructed by the well-known Bogomolov-Tian-Todorov theorem~\cite{Bogomolov,Tian,Todorov}.

For the proof of Theorem~\ref{deform-pK} we first construct a complex nilmanifold $X$ of complex dimension~$4$ with (non-flat) neutral Calabi-Yau structures (see Proposition~\ref{neutralCY-nil}). An interesting class of neutral Calabi-Yau nilmanifolds, called Kodaira manifolds, was constructed by Fino, Pedersen, Poon and S{\o}rensen in \cite{FPPS}. It is worth to note that the complex structure of the nilmanifold $X$ is of a very different kind. Indeed, $X$ has the special feature that it is 4-step and the center of its underlying Lie algebra has minimal dimension, which implies that $X$ is far from being the total space of a torus bundle over a torus. Moreover, in Proposition~\ref{pK-en-SnN} we show that
$X$ provides counterexamples to a conjecture in~\cite{CFU} about the type 
of invariant complex structures on nilmanifolds that admit compatible pseudo-K\"ahler metrics. 
By appropriately deforming the complex nilmanifold $X$, we construct 
a holomorphic family of compact complex manifolds
$\{X_\mathbf{t}\}_{\mathbf{t}\in\Delta}$, with $X_0=X$, 
showing that the neutral Calabi-Yau and 
neutral K\"ahler properties 
are unstable.

In Section~\ref{pseudo-HS} we consider pseudo-Hermitian-symplectic structures, which are an indefinite version of Hermitian-symplectic structures. This geometry naturally arises from small deformations of pseudo-K\"ahler manifolds (see Proposition~\ref{openness-p-HS}). Motivated by a question of Streets and Tian in the positive-definite case~\cite[Question 1.7]{ST}, we prove that there are   
compact complex manifolds admitting pseudo-Hermitian-symplectic structures but 
no pseudo-K\"ahler metrics.

\section{Pseudo-K\"ahler manifolds}\label{general-results}

\noindent
This section starts showing that compact pseudo-K\"ahler manifolds of complex dimension $n\geq 3$ are not stable under small holomorphic deformations of the
complex structure. This motivates the study of conditions under which a sufficiently small deformation $X_{\mathbf t}$ of a pseudo-K\"ahler manifold $X$ is again pseudo-K\"ahler. Here, 
we present a stability result related to the
Bott-Chern cohomology of complex manifolds.
Several explicit constructions of pseudo-K\"ahler nilmanifolds and solvmanifolds
are provided along the section,
illustrating the behaviour of the pseudo-K\"ahler property under their small holomorphic deformations.
In the final part of the section, 
we focus on the behaviour of 
small deformations of 
pseudo-K\"ahler compact complex surfaces.

\smallskip

Let us recall that a \emph{complex analytic family}, or \emph{holomorphic family}, of compact complex manifolds is a proper holomorphic submersion $\pi\colon {\mathcal X} \longrightarrow \Delta$ between two complex manifolds~${\mathcal X}$ and~$\Delta$~\cite{Kod86}. This implies that the fibres $X_{\mathbf{t}} = \pi^{-1}(\mathbf{t})$ are compact complex manifolds of the same dimension. By a classical result of Ehresmann \cite{Ehr47}, any such family is locally $\mathcal{C}^\infty$ trivial (globally, if $\Delta$ is contractible), so all the fibres $X_{\mathbf{t}}$ have the same underlying $\mathcal{C}^\infty$ manifold $M$. 
Consequently, the holomorphic family can be viewed as a collection $\{X_{\mathbf{t}}\}_{\mathbf{t}\in \Delta}$ of complex manifolds $X_{\mathbf{t}}=(M, J_{\mathbf{t}})$, where $J_{\mathbf{t}}$ is the complex structure of $X_{\mathbf{t}}$ for $\mathbf{t}\in \Delta$.

A classical result of Kodaira and Spencer~\cite{KS60} asserts that  
the property of being a compact K\"ahler manifold is stable
under holomorphic deformations. In the following section we prove that such a stability result cannot be extended to compact pseudo-K\"ahler manifolds.

\subsection{Instability of the pseudo-K\"ahler property}\label{ii-1}

Here we prove that compact pseudo-K\"ahler manifolds of complex dimension $n\geq 3$ are 
in general not stable under small deformations of the
complex structure. 
In the proof we will consider 
an appropriate holomorphic family consisting of complex nilmanifolds.
We recall that a nilmanifold is a compact quotient $N=\nilm$ of a 
connected and
simply connected nilpotent Lie group~$G$
by a lattice $\Gamma$ of maximal rank in~$G$. 
A \emph{complex nilmanifold} $X=(N,J)$ is a nilmanifold $N=\nilm$ endowed with an invariant complex structure~$J$, i.e.~$J$ comes from a left-invariant complex structure on~$G$ by passing to the quotient~$\nilm$. 
(For results on complex nilmanifolds see for instance \cite{Ang-libro,COUV,Salamon} and the references therein.)

\begin{proposition}\label{deform-pK-dim6}
There is a holomorphic family of compact complex manifolds 
$\{X_\mathbf{t}\}_{\mathbf{t}\in\Delta}$ of complex dimension $3$,
where $\Delta=\{ \mathbf{t}\in \mathbb{C}\mid |\mathbf{t}|< 1 \}$, such that:
\begin{enumerate}
\item[\rm (i)] $X_0$ is a pseudo-K\"ahler manifold;
\item[\rm (ii)] $X_\mathbf{t}$ does not admit pseudo-K\"ahler metrics for any
$\mathbf{t}\not= 0$.
\end{enumerate}
Hence, the pseudo-K\"ahler property is not stable
under deformations of the complex structure.
\end{proposition}

\begin{proof}
Let $X$ be the complex nilmanifold defined by the following complex structure equations
\begin{equation}\label{ecccus}
d\omega^1 = d\omega^2=0, \quad \ d\omega^3= \omega^{1\bar{2}}.
\end{equation}
Here $\omega^k$ has bidegree (1,0) and $\omega^{j\bar{k}}$ denotes the (1,1)-form $\omega^{j}\wedge\overline{\omega^{k}}$.

Observe that the compact complex manifold $X$ has pseudo-K\"ahler metrics. For instance,
the 2-form
$$
F= i\,\omega^{1\bar{1}} + \omega^{2\bar{3}} - \omega^{3\bar{2}}
$$
satisfies that $F=\bar F$, i.e. the 2-form $F$ is real, and $F^3 \not= 0$, i.e. it is non-degenerate. The form $F$ is compatible with the complex structure $J$ of $X$ because it is of pure type (1,1) on $X$. 
Moreover, from \eqref{ecccus} we get that $dF=0$, so the form $F$ is symplectic.
Thus, $F$ defines a pseudo-K\"ahler metric on $X$.

Notice that the (0,1)-form $\omega^{\bar{3}}$ is $\delbar$-closed on $X$.
We will use the Dolbeault cohomology class
$[\omega^{\bar{3}}] \in H^{0,1}_{\db}(X)$
to perform an appropriate holomorphic deformation of $X$.
For each $\mathbf{t}\in \mathbb{C}$ such that $|\mathbf{t}|<1$, let us consider the complex nilmanifold $X_{\mathbf{t}}$
defined by the following complex basis of (1,0)-forms:
\begin{equation}\label{rel-tt}
\omega_{\mathbf{t}}^1:=\omega^1,\quad \omega_{\mathbf{t}}^2:=\omega^2,\quad \omega_{\mathbf{t}}^3:=\omega^3 + \mathbf{t}\,\omega^{\bar{3}}.
\end{equation}
It follows from \eqref{ecccus} and \eqref{rel-tt} that the complex structure equations for $X_{\mathbf{t}}$ are
\begin{equation}\label{ecccus-t}
d\omega_{\mathbf{t}}^1=d\omega_{\mathbf{t}}^2=0,\quad \
d\omega_{\mathbf{t}}^3= \omega_{\mathbf{t}}^{1\bar2} - \mathbf{t}\,\omega_{\mathbf{t}}^{2\bar1}.
\end{equation}

In order to prove that $X_{\mathbf{t}}$ has no pseudo-K\"ahler metrics for any $\mathbf{t}\not= 0$, we use
Nomizu's theorem \cite{Nomizu} for the de Rham cohomology of nilmanifolds together with the symmetrization
process introduced in \cite{Bel}. More concretely, we take into account that, by \cite[Remark 5]{U}, any closed $k$-form
$\alpha$ on a nilmanifold is cohomologous to the invariant $k$-form $\widetilde\alpha$
obtained by the symmetrization process. 

Fix any ${\mathbf{t}}\in\Delta$. 
Since the complex structure on $X_{\mathbf{t}}$ is invariant,
if there exists a pseudo-K\"ahler structure~$\Theta$ on $X_{\mathbf{t}}$, then the 
form~$\widetilde\Theta$,
obtained by symmetrization of $\Theta$, would be an invariant closed real 2-form of bidegree (1,1)
such that $[\widetilde\Theta]=[\Theta]$  in the second de Rham cohomology group $H_{\rm dR}^2(X_{\mathbf{t}};\mathbb{R})$. 
In particular, this would imply $[{\widetilde\Theta}^3]=[\Theta]^3 \not=0$ in the de Rham cohomology group $H_{\rm dR}^6(X_{\mathbf{t}};\mathbb{R})$, due to the non-degeneracy of $\Theta$ on 
the compact complex manifold
$X_{\mathbf{t}}$.
However, for $\mathbf{t}\not= 0$, it follows from \eqref{ecccus-t} that any invariant closed (1,1)-form $\zeta$ on $X_{\mathbf{t}}$ 
satisfies
$$
\zeta\ \in\  \C\langle \omega_{\mathbf{t}}^{1\bar1},\ \omega_{\mathbf{t}}^{1\bar2},\ \omega_{\mathbf{t}}^{2\bar1},\ \omega_{\mathbf{t}}^{2\bar2} \rangle,
$$
which implies that $\zeta^3=\zeta\wedge\zeta\wedge\zeta=0$, i.e. $\zeta$ is degenerate.
In conclusion, there are no (invariant or not) pseudo-K\"ahler metrics on $X_{\mathbf{t}}$ for $\mathbf{t}\not= 0$.
\end{proof}

\begin{remark}\label{remark-ii-1}
{\rm 
Let us consider $Y_\mathbf{t}=X_\mathbf{t}\times {\mathbb T}^m$, where 
$\{X_\mathbf{t}\}_{\mathbf{t}\in\Delta}$ 
is the holomorphic family given in Proposition~\ref{deform-pK-dim6} 
and ${\mathbb T}^m$ denotes the $m$-dimensional complex torus endowed with its standard K\"ahler metric.
A similar argument as the one in the proof of the previous proposition shows that $Y_\mathbf{t}$ does not admit any pseudo-K\"ahler metric for
$\mathbf{t}\not= 0$. Hence, the pseudo-K\"ahler property is not stable in any complex dimension $n\geq 3$
(see Theorem~\ref{deform-pK} for an irreducible example in complex dimension 4).
In contrast, we will prove in Section~\ref{surfaces} that compact pseudo-K\"ahler surfaces are stable. 
}
\end{remark}

The following example illustrates that, although the pseudo-K\"ahler property is in general not stable, one can find certain deformations
where the existence of pseudo-K\"ahler metrics is preserved.

\begin{example}\label{ejemplo2-full}
{\rm
Let us consider the differentiable family $\{X_t\}_{t\in (-1,1)}$ of compact complex nilmanifolds
defined by the complex structure equations
$$
d\omega_t^1=0,\quad d\omega_t^2=\omega_t^{1\bar{1}},\quad d\omega_t^3= \omega_t^{12} + t\,\omega_t^{1\bar2}.
$$
We note that this is a differentiable family of deformations of the compact complex nilmanifold $X=X_0$ determined by the equations
\begin{equation}\label{X-gen-ecus}
d\omega^1=0,\quad d\omega^2=\omega^{1\bar{1}},\quad d\omega^3= \omega^{12}.
\end{equation}
The manifolds $X_t$ in this family are all pseudo-K\"ahler since, for instance, 
$F_t=i\,(\omega_t^{1\bar{3}}+\omega_t^{3\bar{1}}) +i\,(1+t)\,\omega_t^{2\bar{2}}$ 
defines
a pseudo-K\"ahler metric on $X_t$ for every $t\in (-1,1)$.
}
\end{example}

Hence, one would like to study additional conditions under which the pseudo-K\"ahler property becomes stable. We next establish a condition in terms of the Bott-Chern cohomology.

\subsection{Bott-Chern cohomology and stability of pseudo-K\"ahler manifolds}\label{ii-5}

Here, we show that the stability of the pseudo-K\"ahler property 
is closely related to the variation of the Bott-Chern cohomology.
We recall that the Bott-Chern and the Aeppli cohomologies \cite{Aeppli, Bott-Chern} 
(see also~\cite{Ang-libro}) of a compact complex manifold $X$
are defined, respectively, by
$$
H^{\bullet,\bullet}_{\rm BC}(X):=\frac{\ker\del\cap\ker\delbar}{\imm\del\delbar}
\qquad
\text{ and }
\qquad
H^{\bullet,\bullet}_{\rm A}(X):=\frac{\ker\del\delbar}{\imm\del+\imm\delbar} .
$$
The dimensions of these cohomology groups will be denoted by $h^{p,q}_{\rm BC}(X) = \dim_{\C} H^{p,q}_{\rm BC}(X)$ and $h^{p,q}_{\rm A}(X) = \dim_{\C} H^{p,q}_{\rm A}(X)$.

Suppose that $X$ admits a pseudo-K\"ahler metric 
defined by $F$.
Since the real form $F^k$ is closed and has bidegree $(k,k)$, it defines a Bott-Chern class $[F^k] \in H^{k,k}_{\rm BC}(X)$ for any $1 \leq k \leq n$.
Moreover:

\begin{lemma}\label{necessary-BC-pK-cond}
Let $X$ be a compact complex manifold with $\dim_{\mathbb C} X=n$.
If $X$ admits a pseudo-K\"ahler metric, then
$h^{k,k}_{\rm BC}(X)\geq 1$
for any $1 \leq k \leq n$.
\end{lemma}

\begin{proof}
Let us consider a pseudo-K\"ahler metric on $X$ defined by $F$
and suppose that the class $[F^k] \in H^{k,k}_{\rm BC}(X)$ is zero for some $1 \leq k \leq n$.
This fact implies that $F^k=\del\delbar \beta$ for some form $\beta \in \Omega^{k-1,k-1} (X)$, so
$$
F^n=F^{n-k} \wedge F^k = F^{n-k} \wedge \del\delbar \beta = \del (F^{n-k} \wedge \delbar \beta)
= d (F^{n-k} \wedge \delbar \beta),
$$
which contradits the non-degeneracy of the form $F$.
\end{proof}

The following proposition is an extension 
of the Kodaira-Spencer stability result \cite[Theorem~15]{KS60} for K\"ahler metrics.
We recall that, for every $(p,q)$, the function $t\mapsto h^{p,q}_{\rm BC}(X_t)$ is upper-semi-continuous~\cite{schweitzer}.

\begin{proposition}\label{BC-constant}
Let $X$ be a compact pseudo-K\"ahler manifold, and let $\{X_t\}_{t\in(-\varepsilon,\varepsilon)}$ be a differentiable family of deformations of $X=X_0$,
where $\varepsilon>0$. Suppose that the 
upper-semi-continuous 
function $t\mapsto h^{1,1}_{\rm BC}(X_t)$ is constant.
Then, the compact complex manifold $X_t$ admits a pseudo-K\"ahler metric for any $t$ close enough to~$0$.
\end{proposition}

\begin{proof}
Let $\{\omega_t\}_{t}$ be a family of Hermitian metrics on $X_t$ for $t\in(-\varepsilon,\varepsilon)$.
For each $t$, we consider the Bott-Chern Laplacian $\Delta^{\rm BC}_t$ associated to the Hermitian metric $\omega_t$ on $X_t$ and the corresponding Green operator~$G_t$~\cite{schweitzer}. Denote by $H_t\colon \Omega^{*}_{\C}(X_t) \to \ker\Delta^{\rm BC}_t$ the projection onto the space of harmonic forms with respect to $\Delta^{\rm BC}_t$ (and with respect to the 
$L^2_{\omega_t}$-orthogonal decomposition induced by $\omega_t$~\cite{schweitzer}), 
and by $\pi^{1,1}_t \colon \Omega^{*}_{\C}(X) \to \Omega^{1,1}(X_t)$
the projection of the space of complex forms on $X$ onto the space of $(1,1)$-forms on $X_t$. 

Now, for any $t\in(-\varepsilon,\varepsilon)$, the operator $\Pi_t$ defined by
$$
\Pi_t \;:=\; \left( H_t + \del_t\delbar_t\left(\del_t\delbar_t\right)^{*_t}G_t \right) \circ \pi^{1,1}_t \colon \Omega^{*}_{\C}(X)  \longrightarrow \ker\del_t\cap\ker\delbar_t \;,
$$
gives the projection of the space of complex forms on $X$ onto the space of $\del_t$-closed and $\delbar_t$-closed $(1,1)$-forms on the compact complex manifold $X_t$. Here, $*_t$ is the Hodge-operator with respect to the Hermitian metric $\omega_t$ on $X_t$.
Since the function $t\mapsto h^{1,1}_{\rm BC}(X_t)$ is constant, by elliptic theory \cite[Theorem 7]{KS60} one has that the family $\{\Pi_t\}_{t}$ is smooth in $t$. Let $F_0$ be a pseudo-K\"ahler metric on $X_0=X$.
For $t\in(-\varepsilon,\varepsilon)$, we set
$$
F_t \;:= \frac{\Pi_t F_0 + \overline{\Pi_t F_0}}{2}.
$$
Then, the family $\{F_t\}_t$ is smooth in $t$,
and each $F_t$ is a real (1,1)-form on $X_t$ which is $d$-closed, because it is closed by $\del_t$ and $\delbar_t$.
Since $F_0^n$ is a nowhere vanishing $(n,n)$-form, we have that $F_t$ is non-degenerate for $t$ close enough to $0$.
Therefore, the form $F_t$ defines a pseudo-K\"ahler metric on the compact complex manifold $X_t$ for any $t$ close enough to $0$.
\end{proof}

\begin{example}\label{ejemplo1-BC}
{\rm
We here show that for the holomorphic family 
$\{X_{\mathbf{t}}\}_{\mathbf{t}\in\Delta}$ 
constructed in the proof of Proposition~\ref{deform-pK-dim6},
the upper-semi-continuous function $\mathbf{t}\mapsto h^{1,1}_{\rm BC}(X_{\mathbf{t}})$ is not constant.
According to Proposition~\ref{BC-constant}, this fact gives a reason for the instability
of the pseudo-K\"ahler property along the deformation 
of $X_0$.
Indeed,
for the compact complex manifold $X_0$ we have
$$
H^{1,1}_{\rm BC}(X_0)= \langle \,[\,i\omega^{1\bar{1}}\,],\,[\,i\omega^{2\bar{2}}\,],\,
[\,\omega^{1\bar{2}}-\omega^{2\bar{1}}\,],\,[\,i\,(\omega^{1\bar{2}}+\omega^{2\bar{1}})\,],\,
[\,\omega^{2\bar{3}}-\omega^{3\bar{2}}\,],\,
[\,i\,(\omega^{2\bar{3}}+\omega^{3\bar{2}})\,] \, \rangle,
$$
and,
for any ${\mathbf{t}} \in \Delta-\{0\}$, the Bott-Chern cohomology group of bidegree (1,1) of $X_{\mathbf{t}}$ is
$$
H^{1,1}_{\rm BC}(X_{\mathbf{t}}) = \langle \,[\,i\omega_{\mathbf{t}}^{1\bar{1}}\,],\,[\,i\omega_{\mathbf{t}}^{2\bar{2}}\,],\,
[\,\omega_{\mathbf{t}}^{1\bar{2}}-\omega_{\mathbf{t}}^{2\bar{1}}\,],\,[\,i\,(\omega_{\mathbf{t}}^{1\bar{2}}+\omega_{\mathbf{t}}^{2\bar{1}})\,] \, \rangle.
$$
Therefore, $h^{1,1}_{\rm BC}(X_0)=6$, and $h^{1,1}_{\rm BC}(X_{\mathbf{t}})=4$, for any ${\mathbf{t}} \not=0$.
}
\end{example}

\begin{example}\label{ejemplo2-BC}
{\rm
Let us consider the differentiable family $\{X_t\}_{t\in (-1,1)}$ of compact pseudo-K\"ahler manifolds given in Example~\ref{ejemplo2-full}. 
The Bott-Chern cohomology group of bidegree (1,1) of~$X_t$ is
$$
H^{1,1}_{\rm BC}(X_{t}) = \langle \,[\,i\,\omega_{t}^{1\bar{1}}\,],\, [\,\omega_{t}^{1\bar{2}}-\omega_{t}^{2\bar{1}}\,],\,[\,i\,(\omega_{t}^{1\bar{2}}+\omega_{t}^{2\bar{1}})\,],\,
[\,i\,(\omega_t^{1\bar{3}}+\omega_t^{3\bar{1}}) +i\,(1+t)\,\omega_t^{2\bar{2}}\,]\,\rangle.
$$
Therefore, $h^{1,1}_{\rm BC}(X_{t})=4$ for every $t$, i.e. the function $t\mapsto h^{1,1}_{\rm BC}(X_t)$ is constant. 
}
\end{example}

We next show that the result in Example~\ref{ejemplo2-BC} can indeed be extended to any family 
of deformations of~$X=X_0$.

\begin{proposition}\label{example-general}
Let $X$ be the 
compact pseudo-K\"ahler manifold in Example~$\ref{ejemplo2-full}$ defined by \eqref{X-gen-ecus},
and 
let $\{X_t\}_{t\in(-\varepsilon,\varepsilon)}$ be any differentiable family of deformations of $X=X_0$, where $\epsilon>0$.
Then, $X_t$ admits a pseudo-K\"ahler metric for every $t$ close enough to $0$.
\end{proposition}

\begin{proof}
Since $X$ is a nilmanifold endowed with an invariant complex structure, 
one has by~\cite[Theorem 2.6]{rollenske} that the complex structure $J_t$ 
of any sufficiently small deformation $X_t$ of $X$ is also invariant. 
The dimension of the Bott-Chern cohomology groups of any invariant complex structure $J_t$ 
is given in \cite[Table 2]{AFR} and \cite[Appendix 6]{LUV}. Note that the Lie algebra underlying the nilmanifold~$X$ is precisely $\frh_{15}$, and $h^{1,1}_{\rm BC}(\frh_{15},J)\geq 4$ for any complex structure $J$ on $\frh_{15}$. 

Since $h^{1,1}_{\rm BC}(X_{t})=h^{1,1}_{\rm BC}(\frh_{15},J_t)\geq 4$ varies upper-semi-continuously along differentiable families~\cite{schweitzer}, for any $t$ close enough to $0$ we have
$$
4=h^{1,1}_{\rm BC}(X_{0}) \geq h^{1,1}_{\rm BC}(X_{t})\geq 4,
$$
and thus $h^{1,1}_{\rm BC}(X_{t})= 4$. 
Hence, the function $t\mapsto h^{1,1}_{\rm BC}(X_t)$ is constant and, by Proposition~\ref{BC-constant},  
$X_t$ admits a pseudo-K\"ahler metric for any $t$ close enough to $0$.
\end{proof}

We observe that the condition on $h_{\text{BC}}^{1,1}$ given in Proposition~\ref{BC-constant}
is sufficient but not necessary. To
illustrate this fact, we next study the existence of pseudo-K\"ahler metrics on the small 
deformations of the well-known Nakamura manifold. This manifold is a 
holomorphically parallelizable solvmanifold, i.e. a compact quotient $X=\Gamma\backslash G$ 
where~$G$ is a simply connected \emph{complex}
solvable Lie group and~$\Gamma$ a lattice in~$G$.

Let $G$ be the simply-connected complex solvable Lie group given by
$$
G =\left\{
\begin{pmatrix}e^{z_1}&0&0 &z_2 \\
0&e^{-z_1}&0& z_3 \\
0&0& 1 & z_1 \\
0& 0& 0& 1 \end{pmatrix} \ \mid\ z_1,z_2,z_3 \in \mathbb{C}\right\},
$$
i.e. $G$ is the semi-direct product $G = \mathbb{C} \ltimes_{\varphi} \mathbb{C}^2$, where
 $(z_2,z_3)$ are the coordinates on $\mathbb{C}^2$ and
$$
\varphi(z_1) =
\begin{pmatrix}e^{z_1}&0\\0&e^{-z_1}\end{pmatrix}, \quad 
z_1 \in \mathbb{C}.
$$
One can consider the symplectic form on the Lie group $G$ defined by
\begin{equation}\label{yamada-F-1}
F= i\,dz_1 \wedge d{\bar z}_1 + dz_2 \wedge d{\bar z}_3 + d{\bar z}_2 \wedge dz_3,
\end{equation}
which clearly has bidegree (1,1) with respect to the complex structure of $G$.
That is to say, $F$ is a pseudo-K\"ahler metric on $G$. Notice that $F$ is not left-invariant on $G$, indeed the forms 
$$
\omega^1=dz_1,\quad\  \omega^2=e^{-z_1}\,dz_2,\quad\  \omega^3=e^{z_1}\,dz_3,
$$
constitute a basis of left-invariant
forms of bidegree (1,0) on the complex Lie group $G$.

Yamada proved in~\cite{Yamada-1} that there is a lattice $\Gamma$ of maximal rank in $G$ such that $F$ descends to a
pseudo-K\"ahler metric on $X=\nilm$. 
Let $A \in {\rm SL}(2,\mathbb{Z})$ be a unimodular matrix with distinct real eigenvalues $\lambda$ and $\lambda^{-1}$, and take $a=\log\lambda$. 
One can consider a lattice $\Gamma$ in $G$ of the form $\Gamma = \Gamma_1 \ltimes_{\varphi} \Gamma_{\mathbb{C}^2}$,
with $\Gamma_1=a\,\mathbb{Z}+ 2\,\pi i\,\mathbb{Z}$ and $\Gamma_{\mathbb{C}^2}$ a lattice in $\mathbb{C}^2$.
The compact complex (solv)manifold $X=\nilm$ is known as the (holomorphically parallelizable) Nakamura manifold~\cite{Nakamura}.
Yamada proved in \cite[Theorem 2.1]{Yamada-1} that the symplectic form $F$ given in \eqref{yamada-F-1} descends to $X$,
providing in this way the first  example of a non-toral compact holomorphically parallelizable pseudo-K\"ahler solvmanifold. 
In fact, in terms of the (1,0)-basis $\{\omega^k\}_{k=1}^3$, the form $F$ expresses as
$$
F = i\,\omega^{1}\wedge\omega^{\bar{1}} + e^{2i\,\Imag z_1}\,\omega^{2}\wedge\omega^{\bar{3}}
+ e^{-2i\,\Imag z_1}\,\omega^{\bar{2}}\wedge\omega^{3},
$$
where the functions $e^{2i\,\Imag z_1}$ and $e^{-2i\,\Imag z_1}$ 
are $\Gamma$-invariant. Therefore, $F$ induces a pseudo-K\"ahler metric on 
the Nakamura manifold $X$.

In \cite{Hasegawa-DGA} Hasegawa extended Yamada's result to any compact holomorphically parallelizable solvmanifold of complex dimension 3,
determining all the
lattices of simply-connected unimodular complex solvable Lie groups. 
Hasegawa proves that such a complex solvmanifold admits a pseudo-K\"ahler metric if and only if the Hodge number $h_{\db}^{0,1}=3$ \cite[Theorem 2]{Hasegawa-DGA}. 
Notice that for the Nakamura manifold $X$ one has 
$H_{\db}^{0,1}(X)=\langle [\omega^{\bar{1}}], [e^{-2i\,\Imag z_1}\,\omega^{\bar{2}}], [e^{2i\,\Imag z_1}\,\omega^{\bar{3}}]  \rangle$.

In the following proposition 
we show that there is a small deformation $X_{\mathbf{t}}$ of $X$ 
admitting pseudo-K\"ahler metrics. Note that by~\cite{AK} the Nakamura manifold
has $h_{\text{BC}}^{1,1}(X)=7$, whereas our small deformation satisfies
$h_{\text{BC}}^{1,1}(X_{\mathbf{t}})=3$ for every $\mathbf{t}\neq 0$.

\begin{proposition}\label{nakamura}
There exists a small deformation $X_{\mathbf{t}}$ of the holomorphically parallelizable Nakamura manifold $X$ admitting pseudo-K\"ahler metrics for any ${\mathbf{t}}$.
\end{proposition}

\begin{proof}
In \cite[Section 4]{AK} Angella and Kasuya studied 
some 
deformations of the Nakamura manifold $X=\nilm$. We
here 
consider the deformation given by ${\mathbf{t}}\, \frac{\partial}{\partial z_1}\otimes d\bar{z}_1 \in H^{0,1}(X; T^{1,0}X)$, which corresponds 
to the \emph{case (1)} in their paper.

Note that this deformation defines a holomorphic family $\{X_{\mathbf{t}}\}_{\mathbf{t} \in \Delta}$, 
where $\Delta=\{\mathbf{t}\in \mathbb{C} \mid \, |\mathbf{t}|< 1 \}$, such that $X_0=X$. 
We have the (1,0)-forms on $X_{\mathbf{t}}$ given by
$$
\omega_{\mathbf{t}}^1=
dz_1-\mathbf{t}\,d\bar{z}_1, \quad
\omega_{\mathbf{t}}^2=
e^{-z_1}dz_2, \quad
\omega_{\mathbf{t}}^3=
e^{z_1}dz_3,
$$
whose differentials satisfy 
$$
\left\{
\begin{array}{llll}
d\omega_{\mathbf{t}}^1 \!\!&\!\!=\!\!&\!\! 0,\\[4pt]
d\omega_{\mathbf{t}}^2 \!\!&\!\!=\!\!&\!\! -\frac{1}{1-|\mathbf{t}|^2}\,\omega_{\mathbf{t}}^{12}
+ \frac{\mathbf{t}}{1-|\mathbf{t}|^2}\,\omega_{\mathbf{t}}^{2\bar{1}},\\[6pt]
d\omega_{\mathbf{t}}^3 \!\!&\!\!=\!\!&\!\! \frac{1}{1-|\mathbf{t}|^2}\,\omega_{\mathbf{t}}^{13}
- \frac{\mathbf{t}}{1-|\mathbf{t}|^2}\,\omega_{\mathbf{t}}^{3\bar{1}}.
\end{array}
\right.
$$
A direct calculation shows that the real 2-form $F_{\mathbf{t}}$ defined on $X_{\mathbf{t}}$ by
$$
F_{\mathbf{t}} = i\,\omega_{\mathbf{t}}^{1\bar{1}} + e^{2i\,\Imag z_1}\,\omega_{\mathbf{t}}^{2\bar{3}}
+ e^{-2i\,\Imag z_1}\,\omega_{\mathbf{t}}^{\bar{2}3}
$$
is closed and non-degenerate. Since $F_{\mathbf{t}}$ has bidegree (1,1) with respect to the complex structure on $X_{\mathbf{t}}$, we get a pseudo-K\"ahler metric on the compact complex manifold $X_{\mathbf{t}}$, for any $\mathbf{t} \in \Delta$.
\end{proof}

In complex dimension $3$, there is 
another well-known complex holomorphically parallelizable solvmanifold, 
namely,  
the Iwasawa (nil)manifold. 
Although this manifold is symplectic,
it is proved in~\cite[Theorem 3.2]{CFU-Lefschetz} 
(see also \cite{Yamada-1}) 
that a holomorphically 
parallelizable nilmanifold is pseudo-K\"ahler if and only if it is a complex torus. Hence, the
Iwasawa manifold does not admit pseudo-K\"ahler metrics. Moreover:

\begin{proposition}\label{iwasawa}
The Iwasawa manifold and its small deformations do not admit any pseudo-K\"ahler metric.
\end{proposition}

\begin{proof}
Recall that the Iwasawa manifold $X$ is the compact complex manifold obtained as a quotient of the $3$-dimensional complex Heisenberg group. 
As we have previously stated, by~\cite[Theorem 3.2]{CFU-Lefschetz} 
the Iwasawa manifold is not pseudo-K\"ahler.
Nakamura studied in \cite{Nakamura} the small deformations of the Iwasawa manifold (see also \cite{Ang-libro} for more details).
It turns out that the complex structure equations of any sufficiently small deformation $X_\mathbf{t}$ of the Iwasawa manifold $X=X_0$
can be written as
$$
\left\{
\begin{array}{llll}
d\omega_{\mathbf{t}}^1 \!\!&\!\!=\!\!&\!\! d\omega_{\mathbf{t}}^2 = 0,\\[4pt]
d\omega_{\mathbf{t}}^3 \!\!&\!\!=\!\!&\!\! \sigma_{12}\,\omega_{\mathbf{t}}^{12}
+ \sigma_{1\bar{1}}\,\omega_{\mathbf{t}}^{1\bar{1}} + \sigma_{1\bar{2}}\,\omega_{\mathbf{t}}^{1\bar{2}}
+ \sigma_{2\bar{1}}\,\omega_{\mathbf{t}}^{2\bar{1}} + \sigma_{2\bar{2}}\,\omega_{\mathbf{t}}^{2\bar{2}},
\end{array}
\right.
$$
where the coefficients
$\sigma_{12},\sigma_{1\bar{1}},\sigma_{1\bar{2}},\sigma_{2\bar{1}},\sigma_{2\bar{2}} \in \mathbb{C}$ only depend on
the parameter $\mathbf{t}$ in the deformation space
$\Delta=\{\mathbf{t}=(t_{11},t_{12},t_{21},t_{22},t_{31},t_{32})\in \mathbb{C}^6 \mid \ |\mathbf{t}|< \varepsilon \}$
for a sufficiently small $\varepsilon>0$.
Let $H^+(X_{\mathbf{t}}) 
\subset H_{\rm dR}^2(X_{\mathbf{t}};\mathbb{R})$ be the subspace 
determined by the second de Rham cohomology classes that can be represented by closed real forms 
of bidegree $(1,1)$ on the compact complex manifold $X_{\mathbf{t}}$. 
As proved in
\cite[Proposition 3.4]{LU-ProcAMS}, this subspace is given by 
$$
H^+(X_{\mathbf{t}})=\langle\,[i\,\omega_{\mathbf{t}}^{1\bar{1}}],\,[i\,\omega_{\mathbf{t}}^{2\bar{2}}],
\,[\omega_{\mathbf{t}}^{1\bar{2}}-\omega_{\mathbf{t}}^{2\bar{1}}],\,
[i\,(\omega_{\mathbf{t}}^{1\bar{2}}+\omega_{\mathbf{t}}^{2\bar{1}})]\,\rangle.
$$
It is then clear that any de Rham cohomology class in $H^+(X_{\mathbf{t}})$ is degenerate, so $X_{\mathbf{t}}$
does not admit any pseudo-K\"ahler metric.
\end{proof}


\subsection{Compact pseudo-K\"ahler surfaces}\label{surfaces}

We start this section applying our Proposition~\ref{BC-constant} to provide 
a stability result for the pseudo-K\"ahler property in terms of
the complex invariant~$\Delta^2(X)$ introduced in~\cite{AT4}.

Let us recall that, for each $k\in\mathbb N$, Angella and Tomassini introduced the complex invariant
\begin{equation}\label{def-Deltas}
\Delta^k(X) \;:=\; \sum_{p+q=k} \left( h^{p,q}_{\rm BC}(X) + h^{p,q}_{A}(X) \right) - 2\, b_k, 
\end{equation}
which is a non-negative integer~\cite[Theorem A]{AT4}.
Furthermore, by \cite[Theorem B]{AT4}, a compact complex manifold $X$ satisfies the $\del\delbar$-Lemma
if and only if $\Delta^k(X)=0$ for any $k$. We are interested in the term $\Delta^2(X)$.

\begin{corollary}\label{cor1}
Let $X$ be a compact pseudo-K\"ahler manifold, and let $\{X_t\}_{t\in(-\varepsilon,\varepsilon)}$ be a differentiable family
of deformations of $X=X_0$, where $\varepsilon>0$. If the upper-semi-continuous
function $t\mapsto \Delta^2(X_t)$ is constant,
then $X_t$ admits a pseudo-K\"ahler metric for any $t$ close enough to~$0$.
\end{corollary}

\begin{proof}
Suppose $\Delta^2(X_t)=c$ for any $t\in(-\varepsilon,\varepsilon)$, where $c$ is a non-negative integer.
Then, expanding~\eqref{def-Deltas} for $k=2$, we have
\begin{eqnarray*}
&& h^{2,0}_{\rm BC}(X_0) + h^{1,1}_{\rm BC}(X_0) + h^{0,2}_{\rm BC}(X_0)
+  h^{2,0}_{A}(X_0) + h^{1,1}_{A}(X_0) + h^{0,2}_{A}(X_0) \\[5pt]
  && = \Delta^2(X_0) + 2\, b_{2} = c + 2\, b_{2} = \Delta^2(X_t) +2\, b_{2} \\[5pt]
  && = h^{2,0}_{\rm BC}(X_t) + h^{1,1}_{\rm BC}(X_t) + h^{0,2}_{\rm BC}(X_t)
+  h^{2,0}_{A}(X_t) + h^{1,1}_{A}(X_t) + h^{0,2}_{A}(X_t)\;.
\end{eqnarray*}
Since the functions $t\mapsto h^{p,q}_{\rm BC}(X_t)$ and $t\mapsto h^{p,q}_{A}(X_t)$ are upper-semi-continuous for any $(p,q)$, 
they all must be constant.
In particular, the function $t\mapsto \dim_\C H^{1,1}_{\rm BC}(X_t)$ is constant, so
Proposition~\ref{BC-constant} implies that
the compact complex manifold $X_t$ admits a pseudo-K\"ahler metric for any $t$ close enough to~$0$.
\end{proof}

The following result comes as a direct consequence.
\begin{corollary}\label{cor2}
Any sufficiently small deformation of a compact pseudo-K\"ahler manifold $X$ satisfying
$\Delta^2(X)=0$ admits a pseudo-K\"ahler metric. In particular, any sufficiently small deformation
of a compact pseudo-K\"ahler $\del\delbar$-manifold 
is pseudo-K\"ahler.
\end{corollary}

Let us recall that a compact complex surface is K\"ahler if and only if its first Betti number $b_1$ is even (see Kodaira's classification of surfaces, \cite{Miy74} and \cite{Siu83}, or \cite{Buc99,Lam99} for a direct proof).
To our knowledge, there is no classification of compact complex surfaces admitting 
a pseudo-K\"ahler metric. 
Petean found in \cite{Petean} some 
obstructions to the existence of \emph{indefinite K\"ahler} metrics on
compact complex surfaces in terms of Seiberg-Witten invariants.
Here ``indefinite K\"ahler'' means that the signature of the metric is (2,2), i.e. the metric is pseudo-K\"ahler but non-K\"ahler
(see \cite[Theorem 4]{Petean} for a list of the possible compact complex surfaces 
that might admit an indefinite K\"ahler metric).

\vskip.1cm

Next, we make use of our previous results and a result by Teleman~\cite{Tel06} to prove the stability of the pseudo-K\"ahler property on compact complex surfaces.

\begin{theorem}\label{stability-p-K-surfaces}
Any sufficiently small deformation of a compact pseudo-K\"ahler surface
admits a pseudo-K\"ahler metric.
\end{theorem}

\begin{proof}
Let $X$ be a compact pseudo-K\"ahler surface. By Teleman's result \cite[Lemma 2.3]{Tel06}, 
the invariant $\Delta^2(X)$ 
is either $0$ or $2$.
Moreover, Teleman also proves that $\Delta^2(X)=0$ if and only if the first Betti number $b_1$ is even, which is equivalent to the existence of a K\"ahler metric on $X$. 
Hence, in the case $\Delta^2(X)=0$, any sufficiently small deformation of $X$ is again K\"ahler by \cite{KS60}, so in particular pseudo-K\"ahler.

Let us now focus on the case $\Delta^2(X)=2$. By \cite{Tel06} this is equivalent to $X$ having odd first Betti number, so 
the condition $\Delta^2(X)=2$ is a topological property. Thus,
$\Delta^2(X_t)=2$ for any differentiable family $\{X_t\}_{t\in(-\varepsilon,\varepsilon)}$
of deformations of the compact pseudo-K\"ahler surface $X=X_0$. 
Since the upper-semi-continuous
function $t\mapsto \Delta^2(X_t)$ is constant, by Corollary~\ref{cor1} 
the compact complex surface 
$X_t$ admits a pseudo-K\"ahler metric for any $t$ close enough to~$0$.
\end{proof}

\vskip.1cm

To finish this section we show that there is only one
compact complex non-K\"ahler surface diffeomorphic to a solvmanifold that admits pseudo-K\"ahler metrics, namely the Kodaira-Thurston manifold.

Hasegawa classified in~\cite{Hasegawa} the compact complex surfaces $X$ that are diffeomorphic to a
4-dimensional solvmanifold. Moreover, he proved that the complex structures on such surfaces are invariant 
(see \cite{Ovando} for a study of 4-dimensional Lie algebras admitting pseudo-K\"ahler metrics). 
In fact, in~\cite[Theorem 1]{Hasegawa} it is shown 
that $X$ must be one of the following surfaces:
complex torus, hyperelliptic surface,
Inoue surface of type $\mathcal{S}_M$, primary
Kodaira surface, secondary Kodaira surface, or Inoue surface of type $\mathcal{S}^{\pm}$.
Only the first two are K\"ahler, whereas the other ones have vanishing second Betti number, 
with the only exception of a primary Kodaira surface. 
It is well-known that the latter admits symplectic forms \cite{Thu}.
Consequently, a compact complex non-K\"ahler surface diffeomorphic to a solvmanifold 
admits a symplectic form if and only if it is a primary Kodaira surface.

We recall that a primary Kodaira surface, which we will denote by $KT$, admits pseudo-K\"ahler metrics. By \cite{Hasegawa}, for any complex structure on $KT$ 
there is a global basis $\{\omega^1, \omega^2\}$ of $(1,0)$-forms satisfying
\begin{equation}\label{KT-ecus}
d \omega^1 = 0, \quad  d \omega^2 = \omega^{1\bar{1}}.
\end{equation}
A real $(1,1)$-form $F$ on $KT$ is closed if and only if 
\begin{equation}\label{KT-pseudo-K}
F= i r\, \omega^{1\bar{1}} + u\, \omega^{1\bar{2}} - \bar{u}\, \omega^{2\bar{1}},
\end{equation}
for some $r\in {\mathbb R}$ and $u\in {\mathbb C}$. 
Since $F^2=-2|u|^2\omega^{12\bar{1}\bar{2}}$, we have that $F$ is non-degenerate if and only if $u\not=0$.
Thus, there are pseudo-K\"ahler metrics on $KT$.

As a consequence of our previous discussion, we have:

\begin{proposition}\label{cor3-deform-c-pK-stable}
Let $X$ be a compact complex non-K\"ahler surface diffeomorphic to a solvmanifold.
If $X$ admits a pseudo-K\"ahler metric, then $X$ is a primary Kodaira surface $KT$.
\end{proposition}

\section{Neutral Calabi-Yau manifolds}\label{neutralCY}

\noindent
In this section we focus our attention on a special kind of pseudo-K\"ahler manifolds, namely, neutral Calabi-Yau manifolds. Moreover, the intermediate 
class constituted by 
neutral K\"ahler manifolds 
is also studied. 
We first prove that compact neutral, K\"ahler or Calabi-Yau, surfaces are stable by small deformations of the complex structure.
In higher dimensions, we construct 
an 8-dimensional nilmanifold endowed with a neutral Calabi-Yau metric that 
allows 
us to prove the instability of the neutral K\"ahler and neutral Calabi-Yau properties in any even complex dimension $n\geq 4$. 
It is worth to note that such nilmanifold also provides a counterexample to a conjecture in~\cite{CFU}.

\vskip.1cm 

We first recall some definitions. Let $X=(M,J)$ be a complex manifold of complex dimension $n=2m$.
Following~\cite{FPPS}, a \emph{neutral K\"ahler} structure on $X$ 
is a neutral metric $g$, i.e. of signature $(2m,2m)$, such that

\medskip

$\bullet$  $g$ is compatible with $J$, i.e. $g(JU,JV)=g(U,V)$ for any vector fields $U,V$ on $M$; and

\medskip

$\bullet$  $J$ is parallel with respect to the Levi-Civita connection $\nabla$ of $g$, i.e. $\nabla J=0$.

\medskip

These conditions imply that the 2-form $F(U,V)=g(U,JV)$ is closed, i.e. a neutral K\"ahler structure is in particular pseudo-K\"ahler.

A neutral K\"ahler structure is said to be \emph{neutral Calabi-Yau} if
there exists a nowhere vanishing form $\Phi$ of bidegree $(2m,0)$ with respect to $J$
satisfying $\nabla \Phi=0$. Neutral Calabi-Yau manifolds are Ricci-flat.

\medskip

\subsection{Stability on compact complex surfaces}\label{neutralCY-surfaces}

In this section we study the stability of the neutral K\"ahler and neutral Calabi-Yau properties in complex dimension 2.

\begin{proposition}\label{stability-neutral-Kahler-surfaces}
Let $X$ be a compact complex non-K\"ahler surface. 
Suppose that $X$ admits a neutral K\"ahler metric. Then, 
any sufficiently small deformation of $X$ also
admits neutral K\"ahler metrics.
\end{proposition}

\begin{proof}
By Theorem~\ref{stability-p-K-surfaces}, any sufficiently small deformation  $X_\mathbf{t}$ of $X$ admits pseudo-K\"ahler metrics. Since $X$ is non-K\"ahler by hypothesis, the first Betti number of $X_\mathbf{t}$ is odd, so a pseudo-K\"ahler metric on $X_\mathbf{t}$ cannot have signature (4,0) or (0,4). Hence, the pseudo-K\"ahler metrics are necessarily neutral.
\end{proof}

Let us now observe the following. On the one hand, Petean proves in \cite[Proposition 5]{Petean} that if a compact complex surface $X$ admits a Ricci-flat neutral K\"ahler metric, then its Kodaira dimension 
and its first Chern class are both zero. Moreover, 
$X$ must be a complex torus, a hyperelliptic surface,
or a primary Kodaira surface \cite[Corollary 2]{Petean}.
On the other hand, any compact complex surface with holomorphically trivial canonical bundle is isomorphic to a
K3 surface, a torus, or a primary Kodaira surface. Hence, as a consequence of these results,
one can ensure that
the only compact complex surfaces that can admit 
neutral Calabi-Yau structures are a complex torus and
a primary Kodaira surface. 
We know by \cite{Hasegawa} that 
these two 
are diffeomorphic to a 4-dimensional nilmanifold 
with an invariant complex structure. It is easy to check that they are both neutral Calabi-Yau. 
In fact, the result is clear for the torus, whereas for a primary Kodaira surface it suffices
to remark that, according to the complex equations \eqref{KT-ecus}, the $(2,0)$-form 
$\Phi=\omega^{12}$ is nowhere vanishing and parallel with respect to the Levi-Civita connection of any neutral K\"ahler metric \eqref{KT-pseudo-K}.
Concerning deformations, it is well-known that the
small deformations of the invariant complex structures on the complex torus or on the primary Kodaira surface are again invariant,
so they admit neutral Calabi-Yau structures.

We sum up our previous discussion in the following result.

\begin{proposition}\label{stability-neutral-CY-surfaces}
Let $X$ be a compact neutral Calabi-Yau surface. Then, $X$ is a complex torus or a primary Kodaira surface. 
Moreover, any sufficiently small deformation of $X$ 
admits neutral Calabi-Yau structures.
\end{proposition}

Our next goal is to prove the instability of the neutral Calabi-Yau and neutral K\"ahler properties 
in complex dimension $n\geq 4$. We begin constructing an 8-dimensional nilmanifold endowed with neutral Calabi-Yau metrics. 
To our knowledge, this neutral Calabi-Yau nilmanifold is new and, in addition, it provides counterexamples to a conjecture on pseudo-K\"ahler nilpotent Lie algebras, as we will shortly see.

\subsection{A neutral Calabi-Yau nilmanifold in eight dimensions}\label{neutralCY-8dim-nil}

Let $M$ be a nilmanifold endowed with an invariant complex structure $J$. Let $n$ be the complex dimension of $X=(M,J)$, and suppose that $F$ is an invariant pseudo-K\"ahler metric on $X$. By~\cite{Salamon}, there always exists a closed (non-zero) invariant form $\Phi$ of bidegree $(n,0)$ with respect to $J$, so $\nabla \Phi=0$ for the Levi-Civita connection~$\nabla$ of the invariant metric.
Therefore, any invariant pseudo-K\"ahler metric~$F$ on a complex nilmanifold~$X$ is Ricci-flat (see \cite{FPS}).

\begin{proposition}\label{ejemplo-dim8}
Let $X=(M,J)$ be the $8$-dimensional nilmanifold $M$ endowed with the invariant complex structure $J$ defined by a basis
of $(1,0)$-forms $\{\omega^k\}_{k=1}^4$ satisfying the structure equations
\begin{equation}\label{eleccion}
\left\{\begin{array}{rcl}
d\omega^1 &\!\!\!=\!\!\!& 0,\\[3pt]
d\omega^2 &\!\!\!=\!\!\!& -i\,\omega^{14} + i\,\omega^{1\bar 4},\\[4pt]
d\omega^3 &\!\!\!=\!\!\!& \omega^{12} + \omega^{1\bar 2} - \omega^{2\bar1},\\[4pt]
d\omega^4 &\!\!\!=\!\!\!& - \omega^{1\bar 3} + \omega^{3\bar 1}.
  \end{array}\right.
\end{equation}
Then, $X$ admits pseudo-K\"ahler metrics. Moreover, any invariant pseudo-K\"ahler metric $F$ on $X$ is given by
\begin{equation}\label{pK-concreta}
F=i\,(r\,\omega^{1\bar 1} + s\,\omega^{4\bar 4}) + u\,(\omega^{1\bar 2} - \omega^{2\bar 1})
      + v\,(\omega^{1\bar 3} - \omega^{3\bar 1}) - s\,(\omega^{2\bar 3} - \omega^{3\bar 2}),
\end{equation}
where $r,s,u,v\in\mathbb R$ and $rs\neq 0$.
\end{proposition}

\begin{proof}
It is easy to check that the equations~\eqref{eleccion} satisfy the Jacobi identity, i.e. $d^2 \omega^k=0$ for $1\leq k\leq 4$. Thus, they define a simply-connected, connected, nilpotent Lie group $G$ of real dimension $8$. Moreover, since the structure constants belong to ${\mathbb Q}[i]$, the well-known Malcev's theorem \cite{Malcev} ensures the existence of a lattice $\Gamma$ of maximal rank in $G$. 
We consider the nilmanifold $M=\nilm$, which is endowed by construction with the complex structure $J$ defined by the $(1,0)$-forms $\{\omega^k\}_{k=1}^4$.

Any invariant real $2$-form $F$ of bidegree $(1,1)$ on $X=(M,J)$ can be written as 
\begin{equation}\label{pK-dim8}
F \ = \ \sum_{k=1}^4 i\,x_{k\bar k}\,\omega^{k\bar k}
   \ +\sum_{1\leq k<l\leq 4}\big( x_{k\bar l}\,\omega^{k\bar l}-\bar x_{k\bar l}\,\omega^{l\bar k} \big),
\end{equation}
where $x_{k\bar k}\in\mathbb R$ and $x_{k\bar l}\in\mathbb C$, for $1\leq k, l\leq 4$. Since we are looking for pseudo-K\"ahler metrics, we will study the condition
$dF=0$ and the non-degeneracy condition $F^4\neq 0$.

By a direct computation using the complex structure equations~\eqref{eleccion} we get
$$
dF=\Theta+\overline\Theta,
$$
where $\Theta=\partial F$ is the complex $3$-form of bidegree $(2,1)$ 
given by
\begin{equation*}
\begin{split}
\Theta =& \  2\,i\,\Imag x_{1\bar3}\,\omega^{12\bar1} + 2\,i\,\Imag x_{2\bar3}\,\omega^{12\bar2}
	- (x_{2\bar4}-i\,x_{3\bar3})\,\omega^{12\bar3} + x_{3\bar4}\,\omega^{12\bar4}
	- x_{1\bar4}\,\omega^{13\bar1} + i\,x_{3\bar3}\,\omega^{13\bar2} \\[3pt]
     & - x_{3\bar4}\,\omega^{13\bar3} + 2\,\Imag x_{1\bar2}\,\omega^{14\bar1}
     	+ (x_{2\bar2} - \bar x_{3\bar4})\,\omega^{14\bar2} - i\,(x_{2\bar3} + x_{4\bar4})\,\omega^{14\bar3}
	- i\,x_{2\bar4}\,\omega^{14\bar4} \\[3pt]
& - (x_{2\bar4} + i\,x_{3\bar3})\,\omega^{23\bar1} + (x_{2\bar2} + \bar x_{3\bar4})\,\omega^{24\bar1} + i\,(x_{4\bar4} + \bar x_{2\bar3})\,\omega^{34\bar1}.
\end{split}
\end{equation*}
Now, the closedness condition $dF=0$ is equivalent to $\Theta =0$. 
It is straightforward to check that the latter is satisfied if and only if 
$$
x_{2\bar2}=x_{3\bar3}=x_{1\bar4}=x_{2\bar4}=x_{3\bar4}=\Imag x_{1\bar2}=\Imag x_{1\bar3} = \Imag x_{2\bar3}=0,
\quad \mbox{ and }
    \quad x_{2\bar3}=-x_{4\bar4}.
$$
Replacing these values in~\eqref{pK-dim8} and denoting $x_{1\bar1}=r$, $x_{4\bar4}=s$, $x_{1\bar2}=u$, $x_{1\bar3}=v$,
with $r,s,u,v\in\mathbb R$, one directly gets the expression~\eqref{pK-concreta}.

We finally need to ensure the non-degeneracy condition for $F$. Using~\eqref{pK-concreta}, it is easy to see that
$$
F^4 = - 24\,r s^3\omega^{1234\bar1\bar2\bar3\bar4}.
$$
Therefore, $F^4\neq 0$ if and only if $rs\neq 0$, as stated in the proposition.
\end{proof}

In the following result we prove that the family of pseudo-K\"ahler metrics \eqref{pK-concreta}
provides neutral Calabi-Yau metrics in eight dimensions that are not flat (although they all are Ricci flat).

\begin{proposition}\label{neutralCY-nil}
The complex nilmanifold $X=(M,J)$ constructed in Proposition~$\ref{ejemplo-dim8}$ has non-flat neutral Calabi-Yau structures.
\end{proposition}

\begin{proof}
Let us take the basis of real $1$-forms 
$\{e^1,\ldots,e^8\}$ on $X$ defined by $e^{2k-1}+i\, e^{2k}=\omega^k$, for $1 \leq k \leq 4$, where $\{\omega^1,\ldots,\omega^4\}$ is the basis of $(1,0)$-forms in Proposition~\ref{ejemplo-dim8}
satisfying the complex structure equations~\eqref{eleccion}.
In terms of this real basis, the complex structure $J$ and the pseudo-K\"ahler metrics $F$ given in \eqref{pK-concreta}
express as
$$
\begin{array}{rl}
& J e^1 = -e^2, \quad J e^3 = -e^4, \quad J e^5 = -e^6, \quad J e^7 = -e^8, \\[8pt]
& F = \, 2r\, e^{12} + 2u\, e^{13} + 2v\, e^{15} + 2u\, e^{24} + 2v\, e^{26} - 2s\, e^{35} - 2s\, e^{46} + 2s\, e^{78},
\end{array}
$$
where $r,s,u,v\in\mathbb R$ with $rs\neq 0$. The pseudo-Riemannian metric $g(x,y)=F(Jx,y)$ is 
then given in terms of this real basis by the matrix
\begin{equation}\label{neutral-g}
(g_{ij})_{i,j}= \left(
\begin{array}{cccccccc}
2r & 0  & 0  & -2u & 0  & -2v & 0 & 0 \\[2pt]
0  & 2r & 2u & 0   & 2v & 0   & 0 & 0 \\[2pt]
0  & 2u & 0  & 0   & 0  & 2s  & 0 & 0 \\[2pt]
-2u & 0 & 0  & 0   & -2s & 0  & 0 & 0 \\[2pt]
0  & 2v & 0 & -2s  & 0  & 0   & 0 & 0 \\[2pt]
-2v & 0 & 2s & 0   & 0  & 0   & 0 & 0 \\[2pt]
0 & 0 & 0 & 0 & 0 & 0 & 2s & 0 \\[2pt]
0 & 0 & 0 & 0 & 0 & 0 & 0 & 2s
\end{array}
\right).
\end{equation}
It is easy to see that there are metrics in the family~\eqref{neutral-g} with neutral signature.
In fact, taking for instance $u=v=0$, one has that $rs<0$ is equivalent to the signature being $(4,4)$.
Since the $(4,0)$-form $\Phi=\omega^{1234}$ is parallel,
we conclude that $X$ has neutral Calabi-Yau structures.

Let us now prove that any pseudo-K\"ahler metric~\eqref{pK-concreta} on $X$ is non-flat.
Since the pseudo-K\"ahler structures are invariant, the (complexified) Koszul formula for the Levi-Civita connection $\nabla$
of the metric $g$ reduces to
$$
2 g(\nabla_UV,W)=g([U,V],W)-g([V,W],U)+g([W,U],V),
$$
for (invariant) complex vector fields $U,V,W$ on the complex nilmanifold $X$.
Let $\{Z_j\}_{j=1}^4$ denote the basis of complex vector fields of bidegree (1,0) dual to the
basis $\{\omega^j\}_{j=1}^4$.
Notice that, by complex conjugation, it suffices to compute
$\nabla_{Z_k}Z_j$ and $\nabla_{\overline{Z}_k}Z_j$ for
$1 \leq j,k \leq 4$.
Furthermore, since $\nabla J=0$, 
one has that $\nabla_UV$ is of bidegree (1,0) whenever $V$ is.
In particular, $\nabla_U Z_j$ has type $(1,0)$ for every $1 \leq j \leq 4$.

Let $R$ be the curvature tensor of the pseudo-K\"ahler metric, i.e. $R$ is
given by
$$
R(U,V,W,T)=g\big(\nabla_U\nabla_VW-\nabla_V\nabla_UW-\nabla_{[U,V]}W, \, T\big),
$$
for $U,V,W,T$ complex vector fields on $X$. 
Taking into account the observation in the previous paragraph, complex conjugation and the symmetries of the curvature tensor, one concludes that 
the metric $g$ is non-flat if and only if
$R(Z_i,\bar{Z}_j,Z_k,\bar{Z}_l)\not=0$ for some $i,j,k,l$.
In our case, we will prove that $R(Z_2,\bar{Z}_2,Z_2,\bar{Z}_2)\not=0$.

A direct calculation shows:
$$
\nabla_{\bar{Z}_2}Z_1 = Z_3, \quad\quad \nabla_{\bar{Z}_2}Z_2 = 0, \quad\quad \nabla_{\bar{Z}_2}Z_3 = 0,
\quad
\mbox{ and }
\quad
\nabla_{Z_2}Z_2 = -\frac{is}{r}\, Z_1 -\frac{iv}{r}\, Z_2 +\frac{iu}{r}\, Z_3.
$$
From the complex equations~\eqref{eleccion} we have $[Z_2,\bar{Z}_2]=0$.  
Hence,
$$
\begin{array}{rl}
R(Z_2,\bar{Z}_2,Z_2,\bar{Z}_2)
\!\!&\!\!= \, g(\nabla_{Z_2}\nabla_{\bar{Z}_2}Z_2-\nabla_{\bar{Z}_2}\nabla_{Z_2}Z_2-\nabla_{[Z_2,\bar{Z}_2]}Z_2, \, \bar{Z}_2) \\[5pt]
&= -g(\nabla_{\bar{Z}_2}\nabla_{Z_2}Z_2, \, \bar{Z}_2) \\[5pt]
&=\, g\left(\frac{is}{r} Z_3, \, \bar{Z}_2\right) \\[5pt]
&= - \frac{s^2}{r} \not= 0.
\end{array}
$$
\vskip-.4cm
\end{proof}

Note that in~\cite{FPPS} neutral Calabi-Yau structures on a specific 
class of nilmanifolds
are constructed. This class is given by the so-called 
\emph{Kodaira manifolds}, which are $4m$-dimensional $2$-step nilmanifolds
whose underlying Lie algebras have center of dimension $2m$. 
Moreover, their complex structure
is invariant and preserves the center. Kodaira manifolds are a generalization of the
Kodaira-Thurston manifold $KT$, and they have the structure of a 
principal torus bundle over a torus, with fiber the central torus. 

It is worth to remark that, for the neutral Calabi-Yau nilmanifold $X=(M,J)$ constructed in Proposition~\ref{neutralCY-nil},
the center of the Lie algebra $\frg$ underlying $M$ is 1-dimensional. Hence, it is not invariant under the complex structure $J$ (see 
Section~\ref{counterexample} for more details). 
Furthermore, $\frg$ has nilpotency step equal to~$4$. By \cite{PS}, this implies that
the neutral Calabi-Yau nilmanifold $X=(M,J)$ is far from being the total space of a principal torus bundle over a torus.

\subsubsection{Counterexamples to a conjecture on pseudo-K\"ahler nilmanifolds}\label{counterexample}

\noindent Here we show that the new pseudo-K\"ahler nilmanifold constructed in Proposition~\ref{ejemplo-dim8} provides counterexamples to a conjecture in~\cite{CFU}. 
The conjecture states that an invariant complex structure~$J$ on a nilmanifold $M=\nilm$
must satisfy a certain property so that $(M,J)$ admits pseudo-K\"ahler metrics. 
Let us first recall some results on complex structures on nilpotent Lie algebras (which can be found 
in~\cite{LUV-SnN} and the references therein) and then formulate the conjecture in precise terms.

Let $\frg$ be a nilpotent Lie algebra. 
Complex structures
on $\frg$ 
can be classified into different types attending to
the behaviour of the \emph{ascending $J$-compatible series} of~$\frg$,
which is defined inductively as
$$
\fra_0(J)=\{0\}, \quad 
\quad
\fra_k(J)=\big\{X\in\frg \mid [X,\frg]\subseteq \fra_{k-1}(J)\ {\rm and\ } [JX,\frg]\subseteq \fra_{k-1}(J)\big\}, \ \text{ for } k\geq 1.
$$

Note that $\fra_k(J)$ is an even-dimensional $J$-invariant ideal of $\frg$.
In particular,
$\fra_1(J)$ is the largest subspace
of the center of $\frg$ which is $J$-invariant.

Unlike the usual ascending central series $\{\frg_k\}_k$ of $\frg$, the series $\{\fra_k(J)\}_k$ is adapted to the complex
structure $J$, and it allows to introduce the following partition of the space of complex structures:


\begin{definition}\label{tipos_J}\cite{LUV-SnN}
A complex structure $J$ on a nilpotent Lie algebra $\frg$ is said to be
\begin{itemize}
\item[(i)] \emph{nilpotent}, if there exists an integer~$t>0$ such that~$\fra_t(J)=\frg$;

\item[(ii)] \emph{non-nilpotent}, if $\fra_t(J)\neq \frg$ for every $t\geq 0$; moreover, $J$ is called

\begin{itemize}
 \item[(ii.1)]  \emph{strongly non-nilpotent}, 
if $\fra_1(J)=\{0\}$ $($which implies $\fra_t(J)=\{0\}$ for every~$t)$;
 \item[(ii.2)] \emph{weakly non-nilpotent}, if there is an integer~$t>0$ satisfying~$\{0\}\neq \fra_t(J)=\fra_l(J) \neq \frg$, for every~$l\geq t$.
 \end{itemize}
\end{itemize}
\end{definition}

Notice that $\fra_1(J)\neq \{0\}$ for any nilpotent or weakly non-nilpotent 
complex structure $J$. This allows to construct such structures from other complex
structures defined on lower dimensional nilpotent Lie algebras
(see \cite{LUV-SnN} for details). 
This fact leaves strongly non-nilpotent complex structures as
the essentially new complex structures that arise in each even real dimension. 
In \cite[Section 3.1]{LUV-SnN} it is proved that if $\frg$ admits a 
strongly non-nilpotent complex structure $J$, then the nilpotency step of $\frg$ is at least $3$
(see \cite{LUV-SnN} for other general properties on Lie algebras with strongly non-nilpotent complex structures and structure results up to real dimension 8). 

\vskip.15cm

We can now formulate the following conjecture proposed in \cite{CFU}:

\vskip.15cm

\noindent\textbf{Conjecture \cite[page 123]{CFU}}:  \emph{a complex
structure on a $2n$-dimensional nilpotent Lie algebra must be of nilpotent type in the presence of a
compatible symplectic form}. 

\vskip.15cm

It is proved in \cite{CFU} that the conjecture holds for $n\leq 3$. 
However, the complex structure $J$ given in Proposition~\ref{ejemplo-dim8}
provides a \emph{counterexample} for $n=4$. 
In fact, from the equations \eqref{eleccion} it is straightforward to 
prove that $\fra_1(J)=\{0\}$, that is to say, $J$ is strongly non-nilpotent 
according to Definition~\ref{tipos_J}, and 
it admits the compatible symplectic forms given in \eqref{pK-concreta}. 
Furthermore, 
in every complex dimension $n \geq 4$ we have the following result.

\begin{proposition}\label{pK-en-SnN}
For each $n \geq 4$, there exists a $2n$-dimensional nilmanifold endowed with a non-nilpotent complex structure
that admits pseudo-K\"ahler metrics.
\end{proposition}

\begin{proof}
Let us consider $Y=X\times {\mathbb T}^k$, where $X=(M,J)$ is the $8$-dimensional pseudo-K\"ahler nilmanifold given in Proposition~\ref{ejemplo-dim8} and  ${\mathbb T}^k$ denotes the $k$-dimensional complex torus endowed with any invariant pseudo-K\"ahler metric. Then, $Y$ is a pseudo-K\"ahler nilmanifold of real dimension $2n=8+2k$ with invariant complex structure 
$J_{Y}=J\times J_{{\mathbb T}^k}$ 
satisfying 
$
\{0\}\neq \frb= \fra_t(J_{Y}) \neq \frg\times \frb$, 
for every~$t>0$, where $\frb$ denotes the abelian Lie algebra underlying the torus ${\mathbb T}^k$, and $\frg$ 
the Lie algebra of~$M$. 
In particular, $J_{Y}$ is (weakly) non-nilpotent.
\end{proof}

\subsection{Instability in complex dimension $n\geq 4$}\label{sect4}

\noindent In contrast to the stability results for compact complex surfaces proved in 
Section~\ref{neutralCY-surfaces}, we will next
show that the neutral K\"ahler and neutral Calabi-Yau properties are both unstable in every complex dimension $n\geq 4$.
This constitutes a deep difference with the K\"ahler Calabi-Yau case, for which the deformation space is unobstructed by the well-known Bogomolov-Tian-Todorov theorem \cite{Bogomolov,Tian,Todorov}.
We begin with the following result in complex dimension 4.

\begin{theorem}\label{deform-pK}
There exists a holomorphic family of compact complex manifolds 
$\{X_\mathbf{t}\}_{\mathbf{t}\in\Delta}$ 
of complex dimension $4$, 
where $\Delta=\{ \mathbf{t}\in \mathbb{C}\mid |\mathbf{t}|< 1 \}$, satisfying the following properties:
\begin{enumerate}
\item[\rm (i)] $X_0$ is a neutral Calabi-Yau manifold;
\item[\rm (ii)] $X_\mathbf{t}$ does not admit any pseudo-K\"ahler structure for
$\mathbf{t}\in\Delta\setminus \mathcal{C}$, where $\mathcal{C}$ is the real curve through $0$ given by
$\mathcal{C}=\{ \mathbf{t}\in\Delta \mid \Real\mathbf{t}=0 \}$.
\end{enumerate}
Therefore, neither the neutral Calabi-Yau property nor the neutral K\"ahler property are stable
under deformations of the complex structure.
\end{theorem}

\begin{proof}
The proof is based on an appropriate deformation of the neutral Calabi-Yau 
nilmanifold $X=(M,J)$ found in Proposition~\ref{neutralCY-nil}. 

Let $\{\omega^k\}_{k=1}^4$
be the basis of $(1,0)$-forms on $X$ satisfying \eqref{eleccion}.
Observe that the (0,1)-form $\omega^{\bar{1}}$ defines a Dolbeault cohomology class on $X$. We consider the class
$[\omega^{\bar{1}}] \in H^{0,1}_{\db}(X)$
to perform an appropriate holomorphic deformation of $X$.
For each $\mathbf{t}\in \mathbb{C}$ such that $|\mathbf{t}|<1$, we define the complex structure
$J_\mathbf{t}$ on~$M$ given by the following basis $\{\eta_\mathbf{t}^k\}_{k=1}^4$ of $(1,0)$-forms:
\begin{equation}\label{ec-h-epsilon-deform}
\eta_\mathbf{t}^1:=\omega^1,\quad \eta_\mathbf{t}^2:=\omega^2-\mathbf{t}\,\omega^{\bar1},\quad
\eta_\mathbf{t}^3:=\omega^3,\quad \eta_\mathbf{t}^4:=\omega^4.
\end{equation}
The complex structure equations for $X_\mathbf{t}=(M,J_\mathbf{t})$ are
\begin{equation}\label{ec-h-epsilon-llegada}
\left\{\begin{array}{rcl}
d\eta_\mathbf{t}^1 &\!\!\!=\!\!\!& 0,\\[3pt]
d\eta_\mathbf{t}^2 &\!\!\!=\!\!\!& -i\,\eta_\mathbf{t}^{14}+i\,\eta_\mathbf{t}^{1\bar{4}},\\[4pt]
d\eta_\mathbf{t}^3 &\!\!\!=\!\!\!& \eta_\mathbf{t}^{12}+\mathbf{t}\,\eta_\mathbf{t}^{1\bar1}+\eta_\mathbf{t}^{1\bar2} - \eta_\mathbf{t}^{2\bar1},\\[4pt]
d\eta_\mathbf{t}^4 &\!\!\!=\!\!\!& -\eta_\mathbf{t}^{1\bar3}+\eta_\mathbf{t}^{3\bar1}.
\end{array}\right.
\end{equation}
Observe that the initial structure $J$ is recovered for $\mathbf{t}=0$. Therefore, $X_0=X$ and one immediately
gets part~\textrm{(i)} of the statement. 

To prove~\textrm{(ii)}, since $X_\mathbf{t}$ is a complex nilmanifold, a similar argument as in the proof of Proposition~\ref{deform-pK-dim6} allows us to reduce the study of existence of pseudo-K\"ahler metrics on $X_\mathbf{t}$ to the study of invariant ones.
We first analyze the existence of invariant 
closed $2$-forms $\Omega$ on~$M$ compatible with~$J_\mathbf{t}$, not necessarily real. 
Any such $\Omega$ belongs to $\bigwedge^{1,1}_{J_\mathbf{t}}(\frg^*)$, where $\frg$ denotes the Lie algebra underlying the nilmanifold $M$, and it is given by
$$
\Omega=\sum_{k=1}^4\big(a_k\,\eta_\mathbf{t}^1+b_k\,\eta_\mathbf{t}^2+c_k\,\eta_\mathbf{t}^3+f_k\,\eta_\mathbf{t}^4\big)\wedge\eta_\mathbf{t}^{\bar k},
$$
where $a_k$, $b_k$, $c_k$, $f_k\in\mathbb C$, for $1\leq k\leq 4$.
Making use of the complex structure equations \eqref{ec-h-epsilon-llegada}, one has
$$
d\Omega=\partial_{\mathbf{t}}\Omega + \db_{\mathbf{t}}\Omega,
$$
where 
\begin{equation*}
\begin{split}
\partial_{\mathbf{t}}\Omega &= (a_3-\bar{\mathbf{t}}\,b_3+c_1)\,\eta_\mathbf{t}^{12\bar1} + (b_3+c_2)\,\eta_\mathbf{t}^{12\bar2}
	+ (c_3-b_4)\,\eta_\mathbf{t}^{12\bar3} + c_4\,\eta_\mathbf{t}^{12\bar4}
	- (a_4+\bar{\mathbf{t}}\,c_3)\,\eta_\mathbf{t}^{13\bar1} \\
     &\ \ + c_3\,\eta_\mathbf{t}^{13\bar2} - c_4\,\eta_\mathbf{t}^{13\bar3} - i\,(a_2+b_1-i\,\bar{\mathbf{t}}\,f_3)\,\eta_\mathbf{t}^{14\bar1}
     	+ (f_3-i\,b_2)\,\eta_\mathbf{t}^{14\bar2} - (f_4 + i\,b_3)\,\eta_\mathbf{t}^{14\bar3} \\
     &
	\ \  - i\,b_4\,\eta_\mathbf{t}^{14\bar4} - (b_4+c_3)\,\eta_\mathbf{t}^{23\bar1} - (f_3+i\,b_2)\,\eta_\mathbf{t}^{24\bar1} + (f_4-i\,c_2)\,\eta_\mathbf{t}^{34\bar1}.
\end{split}
\end{equation*}
and
\begin{equation*}
\begin{split}
\db_{\mathbf{t}}\Omega & = - (a_3+c_1-\mathbf{t}\,c_2)\,\eta_\mathbf{t}^{1\bar1\bar2}
	+ (f_1+\mathbf{t}\,c_3)\,\eta_\mathbf{t}^{1\bar1\bar3} -i\,(a_2+b_1+i\,\mathbf{t}\,c_4)\,\eta_\mathbf{t}^{1\bar1\bar4}
     	+ (c_3+f_2)\,\eta_\mathbf{t}^{1\bar2\bar3} \\
     &
	\ \ + (c_4-i\,b_2)\,\eta_\mathbf{t}^{1\bar2\bar4}- (f_4+i\,b_3)\,\eta_\mathbf{t}^{1\bar3\bar4}
     	- (b_3+c_2)\,\eta_\mathbf{t}^{2\bar1\bar2} - c_3\,\eta_\mathbf{t}^{2\bar1\bar3}
	- (c_4+i\,b_2)\,\eta_\mathbf{t}^{2\bar1\bar4} \\
     & \ \ + (f_2-c_3)\,\eta_\mathbf{t}^{3\bar1\bar2}
     	+f_3\,\eta_\mathbf{t}^{3\bar1\bar3} + (f_4-i\,c_2)\,\eta_\mathbf{t}^{3\bar1\bar4}
	- f_3\,\eta_\mathbf{t}^{4\bar1\bar2} - i\,f_2\,\eta_\mathbf{t}^{4\bar1\bar4}.
\end{split}
\end{equation*}
From these expressions we see that $\Omega$ is closed, i.e. $\partial_{\mathbf{t}}\Omega=0=\db_{\mathbf{t}}\Omega$, if and only if
$$
a_4=b_2=b_4=c_3=c_4=f_1=f_2=f_3=0,
$$
and
$$
b_1=-a_2, \quad c_2=-b_3, \quad f_4=-i\,b_3,
\quad -a_3+\bar{\mathbf{t}}\,b_3 =c_1=-a_3-\mathbf{t}\,b_3.
$$
Now, the latter two equalities imply the equation
$$
b_3\,\Real \mathbf{t}=0,
$$
which gives rise to the following two cases.

If $\Real \mathbf{t}=0$, then pseudo-K\"ahler metrics exist on $X_\mathbf{t}$; for instance 
$$
F=i\,\eta_\mathbf{t}^{1\bar 1} - i\,\eta_\mathbf{t}^{4\bar 4} 
-\frac{\mathbf{t}}{2}\,\eta_\mathbf{t}^{1\bar 3} 
+ \eta_\mathbf{t}^{2\bar3}
+ \frac{\bar{\mathbf{t}}}{2}\,\eta_\mathbf{t}^{3\bar 1}
-\eta_\mathbf{t}^{3\bar 2}
$$
is a real closed $(1,1)$-form which is non-degenerate.

However, if we assume $\Real \mathbf{t}\neq 0$, then one has $b_3=0$ and 
consequently $c_2=f_4=0$. Thus, every closed $2$-form $\Omega$ of bidegree (1,1) with respect to the complex structure $J_\mathbf{t}$ is given by
$$
\Omega=a_1\,\eta_\mathbf{t}^{1\bar1} + a_2\,(\eta_\mathbf{t}^{1\bar2}-\eta_\mathbf{t}^{2\bar1}) + a_3\,(\eta_\mathbf{t}^{1\bar3}-\eta_\mathbf{t}^{3\bar1}),
$$
with $a_1, a_2, a_3\in\mathbb C$.

Hence, when $\Real \mathbf{t}\neq 0$, the space of real closed $2$-forms on $\frg$ compatible with $J_\mathbf{t}$, namely 
$\mathcal Z^{+}_{J_\mathbf{t}}(\frg)=\{ \alpha\in\wedge^2 \frg^* \mid d\alpha=0 \text{ and } 
J_\mathbf{t}\alpha=\alpha\} = \{ \alpha\in\wedge^{1,1}_{J_\mathbf{t}}\frg_{\mathbb C}^* \mid d\alpha=0 \text{ and } 
\bar{\alpha}=\alpha\}$,
is generated by
\begin{equation}\label{Re-t-no-nula}
\mathcal Z^{+}_{J_\mathbf{t}}(\frg)=
\langle i\,\eta_\mathbf{t}^{1\bar1},\, \eta_\mathbf{t}^{1\bar2}-\eta_\mathbf{t}^{2\bar1},\, \eta_\mathbf{t}^{1\bar3}-\eta_\mathbf{t}^{3\bar1}\rangle.
\end{equation}
Since
every element in $\mathcal Z^{+}_{J_\mathbf{t}}(\frg)$ is degenerate, no pseudo-K\"ahler metrics
exist on $X_\mathbf{t}=(M, J_\mathbf{t})$ when $\Real \mathbf{t}\neq 0$. 
This clearly implies that $X_\mathbf{t}$ cannot admit any neutral K\"ahler or neutral Calabi-Yau structure for $\mathbf{t}\in\Delta\setminus \mathcal{C}$, where
$\mathcal{C}=\{ \mathbf{t}\in\Delta \mid \Real\mathbf{t}=0 \}$. Hence, both properties are unstable by small deformations.
\end{proof}

In the following theorem we sum up the main results about instability found along the preceding sections.

\begin{theorem}\label{NCY-no-estable}
On compact complex manifolds of complex dimension $\geq 3$, the properties of \emph{``being 
pseudo-K\"ahler''}, \emph{``being neutral K\"ahler''} and \emph{``being neutral Calabi-Yau''}
are not stable under small deformations of the complex structure.
\end{theorem}

\begin{proof}
For the pseudo-K\"ahler property, the result follows from 
Proposition~\ref{deform-pK-dim6} and Remark~\ref{remark-ii-1}.
For the neutral K\"ahler and neutral Calabi-Yau properties in complex dimension 4, the result is given in Theorem~\ref{deform-pK}. In higher dimensions, it suffices to consider the product $Y_\mathbf{t}=X_\mathbf{t}\times {\mathbb T}^{2m}$, where $X_\mathbf{t}=(M, J_\mathbf{t})$ 
is the holomorphic family given in Theorem~\ref{deform-pK} 
and ${\mathbb T}^{2m}$ the $2m$-dimensional complex torus endowed with an invariant neutral Calabi-Yau metric. If $\frg$ and $\frb$ denote the Lie algebras underlying $M$ and ${\mathbb T}^{2m}$, respectively, then the space of 
invariant closed real $2$-forms on $Y_\mathbf{t}$ compatible with the product complex structure $J_{Y_\mathbf{t}}=J_{\mathbf{t}}\times J_{{\mathbb T}^{2m}}$ is given by
$$
\mathcal Z^{+}_{J_{Y_\mathbf{t}}}(\frg\times\frb) = 
\mathcal Z^{+}_{J_\mathbf{t}}(\frg) \oplus 
\left\{ 
\eta_\mathbf{t}^{1} \wedge \bar{\alpha} - \alpha\wedge\eta_\mathbf{t}^{\bar1} \mid \alpha\in\Lambda^{1,0}_{J_{{\mathbb T}^{2m}}}\, \frb^* \right\}
\oplus \mathcal Z^{+}_{J_{{\mathbb T}^{2m}}}(\frb),
$$
where the space $\mathcal Z^{+}_{J_\mathbf{t}}(\frg)$ is described in \eqref{Re-t-no-nula} for any $\mathbf{t} \in \Delta$ such that $\Real \mathbf{t}\neq 0$.
It is easy to see that any element in $\mathcal Z^{+}_{J_{Y_\mathbf{t}}}(\frg\times\frb)$ is degenerate, so 
$Y_\mathbf{t}$ does not admit pseudo-K\"ahler metrics for any
$\mathbf{t}\in\Delta\setminus \mathcal{C}$. 
Since $Y_0=X_0\times {\mathbb T}^{2m}$ is neutral Calabi-Yau, the result follows immediately.
\end{proof}

\section{Pseudo-Hermitian-symplectic structures}\label{pseudo-HS}

\noindent
In this section we consider an indefinite version of the Hermitian-symplectic geometry. The motivation comes from the fact that the small deformations of any pseudo-K\"ahler manifold always posses 
what we will call a pseudo-Hermitian-symplectic structure. 
We will show that there are  
compact complex manifolds with pseudo-Hermitian-symplectic structure but 
not admitting any pseudo-K\"ahler metric.

Recall that a complex structure~$J$ on a symplectic manifold $(M,\Omega)$ is said
to be \emph{tamed} by the symplectic form $\Omega$ if the condition
$\Omega(x,Jx)>0$ is satisfied for all non-zero tangent vectors~$x$. 
Following the terminology of~\cite{ST}, we 
will refer to the pair $(\Omega, J)$ as a \emph{Hermitian-symplectic structure}.
Note that the tamed condition is equivalent to require that the $(1,1)$-component $\Omega^{1,1}$ of the symplectic form $\Omega$ is
positive, i.e. $\Omega^{1,1}$ is a Hermitian metric on the complex manifold $X=(M,J)$. 
No example of a non-K\"ahler compact complex manifold admitting 
a 
Hermitian-symplectic structure is known
(see \cite[page 678]{LZ} and \cite[Question~1.7]{ST}).

By analogy, in the pseudo-Hermitian setting, we introduce the following notion:

\begin{definition}\label{def-p-H-S}
A complex manifold $X=(M,J)$ is called \emph{pseudo-Hermitian-symplectic} if there exists a symplectic form
$\Omega$ on $M$ such that its component of bidegree (1,1) with respect to $J$ is non-degenerate. In such case we will say that the pair $(\Omega, J)$ is a \emph{pseudo-Hermitian-symplectic structure}.
\end{definition}

From the definition, it is clear that any pseudo-K\"ahler manifold is pseudo-Hermitian-symplectic.

\begin{example}\label{ejemplos-KT}
{\rm
For a primary Kodaira surface, with equations given by \eqref{KT-ecus}, 
the pseudo-Hermitian-symplectic (indeed pseudo-K\"ahler) structures are defined in \eqref{KT-pseudo-K}. Observe that the form $\Omega= \omega^{12} + \omega^{\bar{1}\bar{2}}$ is symplectic, but the pair $(\Omega,J)$ is not pseudo-Hermitian-symplectic
because the $(1,1)$-component of $\Omega$ is identically zero.
}
\end{example}

\begin{remark}\label{impli}
{\rm In the positive-definite case, \cite[Proposition 2.1]{EFV} provides an important 
characterization of the Hermitian-symplectic condition. More precisely, 
the existence of such a structure
on a complex manifold  $X$ is shown to be equivalent to the existence of 
a Hermitian metric~$F$ satisfying
$\partial F =\overline \partial\alpha$, for some $\partial$-closed $(2,0)$-form $\alpha$ on $X$. In the following example
we illustrate that a similar result does not hold in the pseudo-Hermitian-symplectic setting.
}\end{remark}

\begin{example}
{\rm
Let us consider a compact complex nilmanifold defined by the complex equations
$$
d\omega^1=d\omega^2=0, \quad d\omega^3= \rho\,\omega^{12}+\omega^{1\bar{1}}+\rho\,\omega^{1\bar{2}}+D\,\omega^{2\bar{2}},
$$
for some $D \in \mathbb{C}$ and $\rho \in \{0,1\}$. 

Let $F=\frac{i}{2} r\,\omega^{1\bar{1}}-\omega^{2\bar{3}}+\omega^{3\bar{2}}$, where $r \in \mathbb{R}^*$.
Then, $F$ is a real form of bidegree $(1,1)$ and
$F^3=3i\, r\, \omega^{1\bar{1}2\bar{2}3\bar{3}} \not=0$.
Moreover, the $(2,0)$-form $\alpha=-\omega^{23}$ is $\partial$-closed
and satisfies
$$
\partial F= \omega^{12\bar{1}} + \rho\,\omega^{12\bar{2}}=-\db \alpha.
$$ 
Take now the real 2-form $\Omega=\alpha+F+\bar{\alpha}$. By the previous condition we have $\partial\alpha=\db \alpha+\partial F=0$,  which implies $d\Omega =0$.

However, the closed 2-form $\Omega$ is not symplectic. In fact,
$$
\Omega^3=(\alpha+F+\bar{\alpha})^3= F^3 + 6 \,\alpha\wedge\bar{\alpha}\wedge F =0.
$$
Even more, if $\rho=0$ and $D\in \mathbb{R}^*$, then the corresponding nilmanifold does not admit any symplectic structure
(indeed, the underlying Lie algebra is $\mathfrak{h}_3$ in the notation of \cite{U}).

For the other values of 
$\rho$ and $D$ 
we have by \cite[Proposition 2.4]{COUV} that the nilmanifold has underlying Lie algebra
$\mathfrak{h}_2$, $\mathfrak{h}_4$, $\mathfrak{h}_6$ or $\mathfrak{h}_8$. These nilmanifolds are symplectic.
}
\end{example}

We next prove that, similarly to the Hermitian case~\cite[Proposition~2.4]{Yang}, the pseudo-Hermitian-symplectic property is open under holomorphic deformations.

\begin{proposition}\label{openness-p-HS}
For compact complex manifolds, the pseudo-Hermitian-symplectic property is stable. 
Therefore, any sufficiently small deformation of a pseudo-K\"ahler manifold is pseudo-Hermitian-symplectic.
\end{proposition}

\begin{proof}
Let us denote by $M$ the real manifold underlying a complex manifold $X$ and by $J$ the complex structure on $M$ such that $X=(M,J)$. 
Let $\Omega$ be a symplectic form on $M$ whose (1,1)-component with respect to $J$ is non-degenerate, i.e. $(\Omega,J)$ is a pseudo-Hermitian-symplectic structure. We consider a holomorphic family of compact complex manifolds $\{X_{\mathbf{t}}=(M, J_{\mathbf{t}})\}_{\mathbf{t}\in \Delta}$, with $\Delta$ containing $0$, such that $X_0=X$.

The symplectic form $\Omega$ decomposes on the compact complex manifold $X_{\mathbf{t}}$ as
$$
\Omega = \alpha_{\mathbf{t}} + F_{\mathbf{t}} + \beta_{\mathbf{t}},
$$
where $\alpha_{\mathbf{t}}$ has bidegree (2,0), $F_{\mathbf{t}}$ is the (1,1) component of $\Omega$, and $\beta_{\mathbf{t}}=\bar\alpha_{\mathbf{t}}$.
By hypothesis, for~$\mathbf{t}=0$ the form~$F_0$ is non-degenerate. 
Hence, one concludes that $F_{\mathbf{t}}$ is also non-degenerate
for any $\mathbf{t} \in \Delta$ sufficiently close to $0 \in \Delta$, so 
$X_{\mathbf{t}}$ satisfies the pseudo-Hermitian-symplectic property.

The second assertion in the proposition is clear since any pseudo-K\"ahler manifold is pseudo-Hermitian-symplectic.
\end{proof}

As we recalled above,  
Streets and Tian pose in \cite[Question 1.7]{ST} the following question, which is still an open problem:  
\emph{Does there exist a compact complex manifold, of complex dimension $\geq 3$, admitting a Hermitian-symplectic structure but no K\"ahler metrics?}

In the following result we prove that the indefinite counterpart of this problem has an  affirmative answer.

\begin{proposition}\label{pseudo-ST}
There exist compact complex manifolds with pseudo-Hermitian-symplectic structure but not admitting any pseudo-K\"ahler metric.
\end{proposition}

\begin{proof}
Let us consider the family of compact complex manifolds
$\{X_\mathbf{t}\}_{\mathbf{t}\in\Delta}$, of complex dimension~$3$, constructed in Proposition~\ref{deform-pK-dim6}. 
We have: 
\begin{enumerate}
\item[\rm (i)] $X_0$ is a pseudo-K\"ahler manifold,
\item[\rm (ii)] $X_\mathbf{t}$ does not admit any pseudo-K\"ahler metric for $\mathbf{t}\not= 0$.
\end{enumerate}
Since $X_0$ is a pseudo-K\"ahler manifold, by Proposition~\ref{openness-p-HS} and (ii) we conclude that, for sufficiently small values of 
$\mathbf{t}\not= 0$, the compact complex manifold $X_\mathbf{t}$ is a pseudo-Hermitian-symplectic manifold with no pseudo-K\"ahler metrics.
\end{proof}

In the next examples we consider pseudo-K\"ahler structures 
on a compact complex manifold $X_0$ and illustrate their behaviour along a small holomorphic deformation $X_{\mathbf{t}}$ of $X_0$.

\begin{example}\label{contraej-p-HS-dim3}
{\rm
Let us consider the family of compact complex manifolds
$\{X_\mathbf{t}\}_{\mathbf{t}\in\Delta=\{ \mathbf{t}\in \mathbb{C} \mid |\mathbf{t}|< 1 \}}$ of complex dimension~$3$ constructed in Proposition~\ref{deform-pK-dim6}.
Using the complex equations \eqref{ecccus} of the complex nilmanifold $X_0$, we have that any invariant pseudo-K\"ahler metric $F$ on $X_0$ is given by
$$
F=i\,(r\,\omega^{1\bar 1} + s\,\omega^{2\bar 2}) + u\,\omega^{1\bar 2} - \bar{u}\,\omega^{2\bar 1}
      + v\,\omega^{2\bar 3} - \bar{v}\,\omega^{3\bar 2},
$$
for some $r,s\in\mathbb R$ and $u,v\in\mathbb C$ satisfying $rv\neq 0$.

A direct calculation using \eqref{rel-tt} shows
that the 2-form $F$ decomposes along the deformation~$X_\mathbf{t}$~as 
$$
F= \alpha_{\mathbf{t}} + F_{\mathbf{t}} + \beta_{\mathbf{t}}
= - \frac{v\,\bar{\mathbf{t}}}{1\!-\!|\mathbf{t}|^2}\,\omega_\mathbf{t}^{23}
+i\,(r\,\omega_\mathbf{t}^{1\bar 1} + s\,\omega_\mathbf{t}^{2\bar 2}) + u\,\omega_\mathbf{t}^{1\bar 2} - \bar{u}\,\omega_\mathbf{t}^{2\bar 1} +
\frac{1}{1\!-\!|\mathbf{t}|^2} (v\,\omega_\mathbf{t}^{2\bar 3} - \bar{v}\,\omega_\mathbf{t}^{3\bar 2})
- \frac{\bar{v}\,\mathbf{t}}{1\!-\!|\mathbf{t}|^2}\,\omega_\mathbf{t}^{\bar 2\bar 3},
$$
where $\alpha_{\mathbf{t}}=- \frac{v\,\bar{\mathbf{t}}}{1-|\mathbf{t}|^2}\,\omega_\mathbf{t}^{23}$  is the (2,0)-component of $F$, and $\beta_{\mathbf{t}}=\bar\alpha_{\mathbf{t}}$. 
The real (1,1)-form $F_{\mathbf{t}}$ is non-degenerate on the manifold $X_\mathbf{t}$, 
i.e. $F$ defines a pseudo-Hermitian-symplectic structure on $X_\mathbf{t}$, which is in accord to Proposition~\ref{openness-p-HS}. 
Notice that $d F_{\mathbf{t}}=\frac{1}{1-|\mathbf{t}|^2} (v\,\bar{\mathbf{t}}\,\omega_\mathbf{t}^{12\bar 2}
+\bar{v}\,\mathbf{t}\,\omega_\mathbf{t}^{2\bar 1\bar 2}) \not=0$ for any $\mathbf{t}\in\Delta\setminus \{0\}$, i.e. the (1,1)-form $F_{\mathbf{t}}$ is not closed. 
Furthermore, by Proposition~\ref{deform-pK-dim6}~(ii), $X_\mathbf{t}$ does not admit any pseudo-K\"ahler structure for
$\mathbf{t}\in\Delta\setminus \{0\}$.
}
\end{example}

\begin{example}\label{contraej-p-HS-dim4}
{\rm
Let $\{X_\mathbf{t}\}_{\mathbf{t}\in\Delta=\{ \mathbf{t}\in \mathbb{C}\mid |\mathbf{t}|< 1 \}}$ be the family of compact complex manifolds of complex dimension~$4$ constructed in Theorem~\ref{deform-pK}.
The manifold $X_0$ is neutral Calabi-Yau, hence pseudo-K\"ahler, and by \eqref{pK-concreta} any invariant pseudo-K\"ahler metric $F$ on $X_0$ is given by
$$
F=i\,(r\,\omega^{1\bar 1} + s\,\omega^{4\bar 4}) + u\,(\omega^{1\bar 2} - \omega^{2\bar 1})
      + v\,(\omega^{1\bar 3} - \omega^{3\bar 1}) - s\,(\omega^{2\bar 3} - \omega^{3\bar 2}),
$$
for some $r,s,u,v\in\mathbb R$ with $rs\neq 0$. A direct calculation using \eqref{ec-h-epsilon-deform} shows that the 2-form $F$ decomposes as follows along the deformation $X_\mathbf{t}$:
$$
F= \alpha_{\mathbf{t}} + F_{\mathbf{t}} + \beta_{\mathbf{t}}
= - s\,\bar{\mathbf{t}}\,\eta_\mathbf{t}^{13}
+i\,(r\,\eta_\mathbf{t}^{1\bar 1} + s\,\eta_\mathbf{t}^{4\bar 4}) + u\,(\eta_\mathbf{t}^{1\bar 2} - \eta_\mathbf{t}^{2\bar 1})
      + v\,(\eta_\mathbf{t}^{1\bar 3} - \eta_\mathbf{t}^{3\bar 1}) - s\,(\eta_\mathbf{t}^{2\bar 3} - \eta_\mathbf{t}^{3\bar 2})
- s\,\mathbf{t}\,\eta_\mathbf{t}^{\bar 1\bar 3},
$$
where $\alpha_{\mathbf{t}}=- s\,\bar{\mathbf{t}}\,\eta_\mathbf{t}^{13}$ is the (2,0)-component of the form $F$, and $\beta_{\mathbf{t}}=\bar\alpha_{\mathbf{t}}$. 
The real (1,1)-form $F_{\mathbf{t}}$ defines a pseudo-Hermitian-symplectic structure on $X_\mathbf{t}$ as it is non-degenerate, accordingly to Proposition~\ref{openness-p-HS}. 
Note that $F_{\mathbf{t}}$ is not closed  for any $\mathbf{t}\in\Delta\setminus \{0\}$, since $d F_{\mathbf{t}}= s (\bar{\mathbf{t}}\,\eta_\mathbf{t}^{12\bar 1}
+\mathbf{t}\,\eta_\mathbf{t}^{1\bar 1\bar 2})$.
Moreover, by Theorem~\ref{deform-pK}~(ii), the compact complex manifold $X_\mathbf{t}$ does not admit any pseudo-K\"ahler structure for every
$\mathbf{t}\in\Delta\setminus \mathcal{C}$, where 
$\mathcal{C}=\{ \mathbf{t}\in\Delta \mid \Real\mathbf{t}=0 \}$.
}
\end{example}

\section*{Acknowledgments}
\noindent
This work has been partially supported by the projects MTM2017-85649-P (AEI/FEDER, UE),
and E22-17R ``Algebra y Geometr\'{\i}a'' (Gobierno de Arag\'on/FEDER).

\end{document}